\documentclass[a4paper]{article}
      \usepackage[english]{babel}
  
      \usepackage[latin1]{inputenc}
      \usepackage{graphicx}
      \usepackage{multicol}
      \usepackage{graphics,epsfig}
      \usepackage{amsmath,amsfonts,amssymb}

      \usepackage{fancyhdr}

      \newtheorem{thm}{Theorem}[section]
      \newtheorem{propo}{Proposition}[section]
      \newtheorem{prop}{Properties}[section]
      \newtheorem{Def}{Definition}[section]
      \newtheorem{rmq}{Remark}[section]
      \newtheorem{lem}{Lemma}[section]
      \newtheorem{nota}{Notation}[section]
      
\title{\bf \large A SEMILINEAR PARABOLIC-ELLIPTIC CHEMOTAXIS SYSTEM WITH CRITICAL MASS 
IN ANY SPACE DIMENSION}
 \author{\normalsize \bf Alexandre MONTARU\\
 \footnotesize Universit\'e Paris 13, Sorbonne Paris Cit\'e, \\
\footnotesize LAGA, CNRS, UMR 7539,\\ 
\footnotesize F-93430, Villetaneuse, France. \\
\small \textit{montaru@math.univ-paris13.fr} }

\date{}
\begin{document}

\maketitle  

\begin{abstract} \footnotesize
We study radial solutions in a ball of $\mathbb{R}^N$  of a semilinear,
parabolic-elliptic Patlak-Keller-Segel system with a nonlinear 
sensitivity involving a critical power. For $N=2$, the latter reduces to the classical "linear"
 model, well-known for its critical mass $8\pi$. We show that a  critical mass 
phenomenon  also occurs  for $N\geq 3$, but with a strongly different qualitative behaviour. More precisely, if the  total mass of cells is  
smaller or equal to the critical mass $\overline{M}$, then the cell density converges to a  
 regular  steady state with support strictly inside the ball as time goes to infinity. In the case of the critical mass, this result is nontrivial 
since
there exists a continuum of stationary solutions and is moreover
in sharp contrast with the case $N=2$ where infinite time blow-up occurs. 
If the total mass of cells is larger than $\overline{M}$, then all radial solutions blow up in finite time. This actually follows from the existence 
(unlike for $N=2$) of a family of self-similar, blowing up solutions with support strictly inside the ball. 
\end{abstract}

\tableofcontents

\section{Introduction}
\subsection{Origin of the problem}

Chemotaxis is the biological phenomenon whereby some cells or bacteria direct 
their movement according to some chemical present in their environment which can 
be attractive or repulsive. 
We shall focus on the case where the chemical is attractive (then called chemoattractant) and self-emitted by cells.
For instance, in case of starvation, amoebas \textit{Dyctyostelium discoideum} emit cyclic adenosine monophosphate (cAMP) which attract themselves. Chemotaxis is thus a strong mean of communication for cells and leads to collective motion. \\For more details on the social life of amoebas \textit{Dyctyostelium discoideum}, see the article \cite{Herrero-Sastre} of M.A. Herrero and L. Sastre. 

\subsubsection{Mathematical formulation}

Assuming that cells and chemoattractant are diffusing and that cells are sensitive to the chemical's concentration gradient (a fact experimentally observed), Patlak in 1953 (cf. \cite{Patlak}) and Keller and Segel in 1970 (cf. \cite{KS}) have proposed the following mathematical model, a parabolic-parabolic system known as Patlak-Keller-Segel system :
\begin{eqnarray}
\rho_t= D_1\,\Delta \rho-\nabla [\chi \, \nabla c]&&\\
c_t=D_2\, \Delta c+\mu \, \rho&&
\end{eqnarray}
where $\rho$ is the cell density, $c$ the chemoattractant concentration, 
$D_1$ and $D_2$ are diffusion coefficients, $\chi=\chi(\rho,c)$ is the sensitivity of cells to the chemoattractant and $\mu$ the creation rate of chemical by cells. \\
 Cells and chemoattractant are assumed to lie in a bounded domain $\Omega$ of $\mathbb{R}^N$ 
with $N\geq 2$ ($\mathbb{R}^2$ or $\mathbb{R}^3$ physically speaking) so we have to specify
 the boundary conditions.\\ For the cell density $\rho$, it is natural to impose a no flux boundary condition
\begin{equation}
D_1\,\frac{\partial \rho}{\partial \nu}-\chi(\rho,c) \, \frac{\partial c}{\partial \nu}=0 \mbox{ \quad  on } \partial \Omega,
\end{equation}
where $\nu$ denotes the outward unit normal vector to the boundary $\partial \Omega$.\\
For the chemoattractant concentration $c$, Dirichlet boundary conditions are assumed :
\begin{equation}
c=0 \mbox{\quad  on }\partial \Omega.
\end{equation}
Some   cells diffuse much slower than the chemoattractant and we will make this 
assumption. 
In this case, two timescales appear
 in the system and to the limit, we can assume that the chemical concentration $c$ 
reaches instantaneously its 
stationary state. After renormalization, these considerations lead to the simplified parabolic-elliptic 
Patlak-Keller-Segel system :
\begin{eqnarray}
\rho_t= \Delta \rho-\nabla [\chi \, \nabla c]&&\\
- \Delta c= \rho&&
\end{eqnarray}
with the same boundary conditions as above, which then become :
\begin{equation}
\frac{\partial \rho}{\partial \nu}-\chi(\rho,c) \, \frac{\partial c}{\partial \nu}=0 \mbox{ \quad  on } \partial \Omega.
\end{equation}
\begin{equation}
c=0 \mbox{\quad  on }\partial \Omega.
\end{equation}

We would like to add that there also exists cells which have a velocity comparable to 
that of the chemoattractant. This is for instance the case of Escherichia coli which 
moreover has  a 'run and tumble' motion. Hence, in this case, 
 diffusion does not seem to be the most
suitable modeling. On that subject, see the article of B. Perthame \cite{Perthame} 
for a kinetic approach
 which takes into account these characteristics and allows to recover 
the Patlak-Keller-Segel model in a diffusion limit.\\

For a review on mathematics of chemotaxis, see the chapter written by M.A. Herrero in \cite{Herrero} and the article \cite{HP} of  T. Hillen and K. J. Painter. For a review on the Patlak-Keller-Segel model, 
see both articles of D. Horstmann \cite{Horstmann1,Horstmann2}.\\

In \cite{Horstmann-Winkler}, D. Horstmann and M. Winkler have studied the case where the 
sensitivity $\chi$ depends only on $\rho$ and shown that :
\begin{itemize}
\item If $\chi(\rho)\leq C \rho^q$ for $\rho\geq 1$ and $q<\frac{2}{N}$, then 
the cell density $\rho$ exists globally and is even uniformly bounded in time.
\item If $\chi(\rho)\geq C \rho^q$ for $\rho\geq 1$ and $q>\frac{2}{N}$, then 
$\rho$ can blow up.
\end{itemize}
See also \cite{DW, LS, MS, Nagai, SK} for related results.\\

Thus, the power $q=\frac{2}{N}$ of the nonlinearity $\chi(\rho)$ is critical for that system. 
This is why we are going to focus on the following problem, noted $(PKS_q)$, 
with a special interest to the case $q=\frac{2}{N}$ :
\begin{equation}
\label{PKS}
(PKS_q)\qquad  \left\lbrace \qquad 
\begin{array}{ccc}
\rho_t= \Delta \rho-\nabla [\rho^q \, \nabla c]&&\\
- \Delta c= \rho&&
\end{array}
\right.
\end{equation}
where the boundary conditions become 
\begin{equation}
\label{BC-rho}
\frac{\partial \rho}{\partial \nu}-\rho^q \, \frac{\partial c}{\partial \nu}=0 \mbox{ \quad  on } \partial \Omega,
\end{equation}
\begin{equation}
\label{BC-c}
c=0 \mbox{\quad  on }\partial \Omega.
\end{equation}

We would like to stress that $q=\frac{2}{N}$ is exactly the exponent for which the mass, i.e. the $L^1$ norm of $\rho$, is invariant by the rescaling of this system, given by
\begin{eqnarray}
\rho_\lambda(t,y)&=&\lambda^{\frac{2}{q}}\;\rho(\lambda^2\,t,\lambda\,y)\\
c_\lambda(t,y)&=&\lambda^{\frac{2}{q}-2}\;c(\lambda^2\,t,\lambda\,y)
\end{eqnarray}
for all $t>0$, $y\in\mathbb{R}^N$ and $\lambda>0$.
This fact opens the door to the possibility of a critical mass phenomenon.

\begin{rmq}
System $(PKS_q)$ can also be seen as the macroscopic description of a collection of $n$ particles 
following a generalized stochastic Langevin equation. Making a mean field approximation, it is actually obtained as a nonlinear Fokker-Planck equation in a proper thermodynamic limit $n\rightarrow \infty$. For more details, see the article of P.H. Chavanis \cite[section 3.4]{Chavanis}.
\end{rmq}

\subsubsection{Radial setting}

In this paper, we restrict our study to the case of radially symmetric solutions of $(PKS_q)$ where $\Omega$ is the open unit ball $B\subset \mathbb{R}^N$ centered at the origin. Note that by using the scaling of the system and its invariance by translation, we can of course cover the case of any open ball of $\mathbb{R}^N$. \\

We would like to point out that for $N=2$, the critical exponent is $q=1$, so the latter system reduces to the most studied Keller-Segel parabolic-elliptic type system :
\begin{eqnarray}
\rho_t= \Delta \rho-\nabla [\rho \, \nabla c]&&\\
- \Delta c= \rho&&
\end{eqnarray}
It is a well-known fact that this system exhibits a critical mass phenomenon. More precisely, 
denoting $\overline{m}$ the total mass of the cells, it has been shown for radially symmetric solutions that :
\begin{itemize}
\item If $\overline{m}<8\pi$, then $\rho(t)$ is global and converges to a steady state as $t$ 
goes to infinity. \\
(see  \cite{BKLN} by P. Biler, G. Karch, P. Lauren\c cot and T. Nadzieja).
\item If $\overline{m}=8\pi$, then $\rho(t)$ blows up in infinite time to a Dirac mass centered at the origin. \\
(see again \cite{BKLN} and  \cite{Kavallaris-Souplet} by N.I. Kavallaris and P. Souplet for refined asymptotics).
\item If $\overline{m}>8\pi$, then $\rho(t)$ blows up in finite time to a Dirac mass.\\
 (see  \cite{Herrero-Velazquez} by M.A. Herrero and J.L. Velazquez).
\end{itemize}
Moreover, this system exhibits a similar phenomenon in the case of the whole plane 
$\mathbb{R}^2$. See the work of \cite{BKLN2, Blanchet2, Blanchet3, Dolbeault}. In the nonradial case in a 
bounded domain, results are slightly different (see the book \cite{Suzuki} of T. Suzuki). The behaviour of the parabolic-parabolic system in $\mathbb{R}^2$ seems more intricate. See \cite{BCD, CALCOR}.\\

From now on, we consider the case $N\geq 3$.\\

Adapting the procedure described in the article \cite{BHN} of P. Biler, D. Hilhorst and T. Nadzieja 
(or also in \cite{BKLN, Kavallaris-Souplet}), we can reduce the system 
$(PKS)_q$ to a single one-dimensional equation. \\
Indeed, denoting $Q(t,r)=\int_{B(0,r)} \rho(t,y)\,dy$ the total mass of the cells in $B(0,r)$ at time $t$ for $0\leq r\leq 1$, we can make the following formal computations :
\begin{align*}
 Q_t&=\int_{S(0,r)}\left[\frac{\partial \rho}{\partial \nu}-
\rho^q\,\frac{\partial c}{\partial \nu}\right]\,d\sigma_N\\
&=\int_{S(0,r)}  \tilde{\rho}_r\,d\sigma_N-\tilde{\rho}^q\,\int_{B(0,r)}\Delta c\,dy\\
&=\sigma_N \,r^{N-1}\tilde{\rho}_r
+\tilde{\rho}^q\,Q
\end{align*}
 where $\sigma_N$ denotes the surface area of the unit sphere in $\mathbb{R}^N$ 
 and  $$\rho(t,y)=\tilde{\rho}(t,|y|)$$ for any $y\in \overline{B}$.
Since we can write $$Q(t,r)=\int_0^r\sigma_N\, s^{N-1}\tilde{\rho}(t,s)ds,$$
we have both following formulas $$Q_r=\sigma_N\,r^{N-1}\,\tilde{\rho}$$
and $$\sigma_N \,r^{N-1}\tilde{\rho}_r=Q_{rr}-\frac{N-1}{r}Q_r,$$
which imply :
\begin{equation}
Q_t=Q_{rr}-\frac{N-1}{r}\,Q_r+\left[\frac{Q_r}{\sigma_N\,r^{N-1}}\right]^q\,Q.
\end{equation}
Then, setting $P(t,x)=Q(t,x^{\frac{1}{N}})$, we obtain :
\begin{equation}
P_t=N^2\, x^{2-\frac{2}{N}}\,P_{xx}+\left[\frac{N}{\sigma_N}\right]^q\,{P_x}^q\,P.
\end{equation}
Finally, setting $u(t,x)=\frac{1}{N^{\frac{2}{q}}\, V_N}P(\frac{t}{N^2},x)$ 
where $V_N=\frac{\sigma_N}{N}$ is the volume of $B$, we get :
\begin{equation}
u_t=x^{2-\frac{2}{N}}\,u_{xx}+u\,{u_x}^q.
\end{equation}
Moreover, by the no flux boundary condition on the cell density, it is formally clear that the total cell mass $\overline{m}$ is constant in time. Hence, setting $$m=\frac{\overline{m}}{N^{\frac{2}{q}}V_N},$$
 we also have  the boundary conditions that for all $t\geq 0$, $$u(t,0)=0$$ 
 $$u(t,1)=m.$$
 A simple calculation also shows that, for $r\geq 0$,  
 \begin{equation}
 \label{relation_rho_derivee_de_u}
 \tilde{\rho}(t,r)=N^{\frac{2}{q}}\,u_x(N^2\,t,r^N).
 \end{equation}
 
Hence, $\rho$ is simply proportional to $u_x$, up to a time rescaling and a change of variable. It means that 
the derivative of $u$ is the quantity with physical meaning and should then be  nonnegative.\\

Finally, we shall focus on the following problem, noted $(PDE_m)$ :
\begin{equation}
\label{equ_u_1}
(PDE_m)\qquad \left\lbrace \qquad 
\begin{array}{ccc}
u_t=x^{2-\frac{2}{N}}\,u_{xx}+u\,{u_x}^q & t>0 & 0<x\leq 1\\
u(t,0)=0 &t\geq 0&\\
u(t,1)=m &t\geq 0&\\
u_x(t,x)\geq 0 &t> 0& 0\leq x\leq 1.
\end{array}
\right.
\end{equation}
Conversely, starting from a solution $u$ of $(PDE_m)$ we would like to show (at least formally) how to get a solution of $(PKS_q)$.\\
First, it is easy to check that 
$$\tilde{\rho}_r(t,r)=N^{\frac{2}{q}+1}r^{N-1}u_{xx}(N^2\,t,r^N)$$
and, denoting $\tilde{c}$ the radial profile of $c$, since $-\Delta c=\rho$, we have 
$$\tilde{c}_r(t,r)=-\frac{N^{\frac{2}{q}-1}}{r^{N-1}}\, u(N^2\,t,r^N)$$
 so that, by (\ref{BC-c}), we obtain 
\begin{equation}
\label{formule_c}
 \tilde{c}(t,r)=\int_r^1 \frac{N^{\frac{2}{q}-1}}{s^{N-1}}\, u(N^2\,t,s^N).
\end{equation} 
Now, we define $(\rho,c)$ by their profiles given in formulas (\ref{relation_rho_derivee_de_u}) and (\ref{formule_c}).\\
If we denote $r=|y|$ for any $y\in \overline{B}$ and 
$$\alpha(r)=\tilde{\rho}_r-\tilde{\rho}^q\,\tilde{c}_r,$$ then, by the following general fact 
$$r^{N-1}\, \mathop{\rm div}[\alpha(r)\vec{e}_r]=\frac{d}{dr}[r^{N-1}\alpha], $$
we obtain 
\begin{align*}
r^{n-1}\, \mathop{\rm div}[\nabla \rho-\rho^q \nabla c]&= \frac{d}{dr}
[r^{N-1}\tilde{\rho}_r-\tilde{\rho}^q r^{N-1}\tilde{c}_r]\\
&=N^{\frac{2}{q}+1} \frac{d}{dr}\left[r^{2N-2}u_{xx}(N^2 t,r^N)+(u\, {u_x}^q)(N^2t,r^N)\right]\\
&=N^{\frac{2}{q}+1} \frac{d}{dr} u_t(N^2t,r^N) \mbox{\qquad  \qquad  \qquad  by }(\ref{equ_u_1})\\
&=r^{N-1}\; N^{\frac{2}{q}+2}  u_{xt}(N^2t,r^N)\\
&=r^{N-1}\rho_t \mbox{ \quad\qquad  \qquad \qquad \qquad \qquad by }(\ref{relation_rho_derivee_de_u}).
\end{align*}
Hence, $(\rho,c)$ is a solution of $(PKS_q)$. We just have to check that $\rho$ also satisfies the no flux boundary conditions (\ref{BC-rho}) which are equivalent to 
\begin{equation}
\label{condition_au_bord_2}
 \tilde{\rho}_r-\tilde{\rho}^q\,\tilde{c}_r=0 \mbox{\qquad for }r=1.
 \end{equation}
Thanks to the previous formulas on $\tilde{\rho}_r$ and $\tilde{c}_r$,  (\ref{condition_au_bord_2}) becomes
$$\left[u_{xx}+{u_x}^q\, u\right](N^2\,t,1)=0, $$
which can be obtained from $(\ref{equ_u_1})$ since $u_t(N^2\,t,1)=0$.\\

Now, it seems  reasonable to consider problem $(PDE_m)$ as our model for chemotaxis.
Equation (\ref{equ_u_1}) presents two difficulties since the diffusion is degenerate at $x=0$ 
and the nonlinearity is not Lipschitz continuous. We shall assume that the initial data $u_0$ 
belongs to the class 
$$Y_m=\{u\in C([0;1]),\:u \mbox{ nondecreasing },\; u'(0) \mbox{ exists, } u(0)=0,\;u(1)=m\}.$$
For such $u_0$, we have established in \cite{Montaru} the existence
 and uniqueness of a maximal classical solution $u$ such that $u(t)\in Y_m$ for all $t\in 
[0,T_{max}(u_0))$, 
where $T_{max}(u_0)$ is the maximal existence time. 
See Subsection \ref{problem_PDE} below for precise definitions and more details.

\subsection{Main results}

\underline{We now focus on the case 
of the critical exponent $q=\frac{2}{N}$ with $N\geq 3$.}\\

The set of stationary solutions can be precisely described.\\ 
We shall prove that the stationary solutions of $(PDE_m)$ are the restrictions to $[0,1]$ 
of a family of functions $(U_a)_{a\geq 0}$ with the following properties :
\begin{itemize}
\item $U_1(0)=0$, $U_1$ is nondecreasing and reaches its maximum $M$ at 
$x=A$ from which $U_1$ is flat.
\item All $(U_a)_{a\geq 0}$ are obtained by dilation of $U_1$, i.e. $U_a(x)=U_1(ax)$ for all $x\geq 0$.
\end{itemize}

Using this, we can then prove :

\begin{thm} Let $m\geq 0$. Considering problem $(PDE_m)$ with $q=\frac{2}{N}$ :
\label{thm_stat}
\begin{itemize}
\item[i)] If $0\leq m < M$, then there exists a unique stationary solution.
\item[ii)] If $m=M$, there exists a continuum of steady states :  $\left(U_a|_{[0,1]}\right)_{a\geq A}$.\\ Note that the corresponding cell densities have their support strictly inside $B$.
\item[iii)] If $m>M$, there is no stationary solution.
\end{itemize}
\end{thm}

The previous theorem leads us to set the following definition.

\begin{Def}\quad \\
We call $M$ the critical mass of problem $(PDE_m)$ for $q=\frac{2}{N}$ and
$$\overline{M}=N^N V_N\times M $$
the critical mass of system $(PKS_q)$.  
\end{Def}

We would like to stress that the system $(PKS_q)$ exhibits two levels of criticality. We have already seen the first level which consists in choosing the right exponent $q=\frac{2}{N}$ in order to balance the diffusion and aggregation forces in the system. 
Once this exponent is chosen, a second level of criticality arises with the choice of the mass. The following two theorems state that a critical mass phenomenon indeed occurs.

\begin{thm}
\label{thm_convergence}
 Let $N\geq 3$ and $q=\frac{2}{N}$.\\
 If $m\leq M$ and $u_0\in Y_m$, then
 $$u(t) \underset{t \rightarrow + \infty}{ \longrightarrow } 
U_a \mbox{ \quad in \quad } C^1([0,1])$$ 
for some $a\geq 0$.
\end{thm}

\begin{thm}
\label{blow-up}
 Let $N\geq3$ and $q=\frac{2}{N}$. \\If $m>M$,  then for all $u_0\in Y_m$, 
$$T_{max}(u_0)<\infty.$$
Moreover, denoting $\mathcal{N}[u]=\underset{x\in(0,1]}{\sup}\frac{u(x)}{x}$, we have :
$$\underset{t\rightarrow T_{max}}{\lim} \mathcal{N}[u(t)]=\infty .$$
\end{thm}

In addition, for slightly supercritical mass, 
we can show the existence of blowing-up self-similar solutions.

\begin{thm}
\label{thm_self-similar}
There exists $M^+>M$ such that for all $m\in (M,M^+]$ there exists a family of blowing-up self-similar solutions of problem $(PDE_m)$. \\
Moreover, the corresponding cell densities have their support strictly inside $B$.
\end{thm}

\subsection{Comments and related results}
\subsubsection{Description of the ideas of the proofs}

The global existence part of Theorem \ref{thm_convergence} for subcritical or 
critical mass is based 
on comparison with suitable supersolutions, combined with some 
continuation results obtained in \cite{Montaru}. Our convergence  statements
 heavily rely on Lyapunov functional type arguments. 

More precisely, we show that the evolution problem $(PDE_m)$ induces 
a gradient type dynamical system on 
$Y_m^1=Y_m\cap C^1([0,1]) $, 
 with global relatively compact trajectories. 
Moreover, we exhibit a strict Lyapunov functional :
$$\mathcal{F}[u]=\int_0^1 \frac{\dot{u}^{2-q}}{(2-q)(1-q)}-\frac{u^2}{2 x^{2-q}} dx.$$ 
Indeed, formally, it is easy to check that
$$\frac{d}{dt} \mathcal{F}[u(t)]=\int_0^1 u_{tx}
\frac{ \dot{u}^{1-q} }{1-q}-\frac{u_t \,u}{x^{2-q}}dx=\left[ u_{t}
\frac{ \dot{u}^{1-q} }{1-q}\right]_0^1-
\int_0^1 u_t\left[\frac{d}{dx}\frac{ \dot{u}^{1-q} }{1-q}+
\frac{ u}{x^{2-q}}\right]dx.$$
Thanks to the boundary conditions and to $(\ref{equ_u_1})$, we then have
$$\frac{d}{dt} \mathcal{F}[u(t)]=-\int_0^1 \frac{(u_t)^2}{\dot{u}^q \, x^{2-q}}dx.$$
However, this computation is not rigorously valid, since $u_x$ can vanish on a 
whole interval for instance. Nevertheless, we can overcome this difficulty
 and prove that $\mathcal{F}$ is indeed a strict Lyapunov functional by expressing it as
 the limit as $\epsilon$ goes to zero of a family of strict Lyapunov functionals 
$\mathcal{F}_\epsilon$ for 
suitable approximate problems (cf. problem $(PDE_m^\epsilon)$ introduced in Subsection 
\ref{problem_PDE_eps}).
We note that the proof of the compactness of the trajectories relies on a 
different transformation, leading to another auxiliary problem $(tPDE_m)$ 
(cf. Subsection \ref{problem_tPDE} below).
In the subcritical case, since there is a single steady state, this
 immediately implies the convergence of the trajectory. But in the 
critical case, the situation is more delicate, since there exists a 
continuum of steady states and the solution could oscillate without 
stabilizing. However, thanks to a good relation between order and 
topology of the set of stationary solutions, we can prove stabilization by arguments in the spirit of (though
simpler than) those in the articles \cite{Matano, Zelenyak} of H. Matano and T.I. Zelenyak.\\

As for our blow-up results (Theorems \ref{blow-up} and \ref{thm_self-similar}), their
 proofs 
rely on the construction
 of a subsolution which becomes a self-similar solution after some time. 
The latter's profile is solution of an appropriate auxiliary 
 ordinary differential equation which is a perturbation of the stationary 
solution's equation. The construction, as well 
as the study of the steady states (cf. Theorem \ref{thm_stat}), requires 
some rather delicate 
ODE arguments.

\subsubsection{Open problems}
\begin{itemize}
 \item[i)] A natural and very interesting question would 
be to determine the basin of attraction of a given steady state $U_a$ with $a \geq A$.
\item[ii)] For the self-similar solutions in Theorem \ref{thm_self-similar}, 
it is easy to see that 
the blow-up 
rate of the central density of cells 
(proportional to $u_x(t,0)$) behaves like $(T_{max}- t)^{-\frac{N}{2}}$. 
It would be interesting to know if all solutions of problem $(PDE_m)$ blow up at the self-similar
rate 
 or if there also exists blow-up of type $II$, i.e. if there exists solutions with blow-up speed faster than that of the self-similar solutions. 
\end{itemize}

\subsubsection{Comparison with the case $N=2$}

It is instructive to compare the cases $N\geq 3$ and $N=2$ for problem $(PDE_m)$ with $q=\frac{2}{N}$. The behaviour 
is the same for both in the subcritical case since solutions converge to a unique steady state and also in the supercritical case since solutions blow up in finite time. But for the critical case, the qualitative behaviour differs strongly. Indeed, for $N=2$, blow-up occurs in infinite time whereas for $N\geq 3$, there is still convergence to a regular steady state whose corresponding cell density has support strictly inside $B$, a phenomenon which never occurs for $N=2$. We would like to suggest an "explanation" for this.\\
Denoting $S_N$ the set of stationary solutions for $N\geq2$, we can see that we could as well define the critical mass as
$$ M=\underset{U\in S_N}{\sup}\|U\|_{\infty,[0,1]}.$$
The main difference is that this supremum is not reached for $N=2$ whereas it is for $N\geq 3$,
 which allows us in the latter case to find a supersolution that prevents blow-up. Thus,
 convergence or infinite-time blow-up seems to be determined by whether or not the critical mass is reached by stationary solutions.

\subsubsection{Related literature for porous medium type diffusion}

Finally, we would like to make the link between our work and the article \cite{BCL} of 
A. Blanchet, J.A. Carrillo and P. Lauren\c cot (see also the articles \cite{Kim_Yao} of I. Kim and Y. Yao and \cite{BRB} of J. Bedrossian, N. Rodriguez and A.L. Bertozzi for further results in this direction). It will allow to identify a 
formula for the critical mass $M$.\\
The authors there study in the whole space $\mathbb{R}^N$ for $N\geq 3$ the following Patlak-Keller-Segel system $(PKS^p)$ with porous-medium like nonlinear diffusion :
\begin{eqnarray}
\mu_t= \mathop{\rm div}[\nabla \mu^p-\mu \nabla c]&t>0&x\in \mathbb{R}^N\\
-\Delta c=\mu &t>0& x\in \mathbb{R}^N
\end{eqnarray}
where $\mu$ is the cell density and $c$ the concentration of the chemoattractant.\\
They could show that for the critical exponent $p=2-\frac{2}{N}$, the system $(PKS^p)$ exhibits a critical mass phenomenon. See also \cite{Cieslak_Winkler} for a  explanation of this exponent for parabolic-elliptic Patlak-Keller-Segel systems with general nonlinear diffusion.\\
More precisely, denoting $\overline{m}$ the total mass of the cells, they have shown the existence of $M_c$ such that :
\begin{itemize}
\item If $\overline{m}<M_c$, solutions exist globally.
\item If $\overline{m}=M_c$, solutions exist globally in time. Moreover, there are infinitely many compactly supported stationary solutions.
\item If $\overline{m}>M_c$, there are solutions which blow up in finite time.
\end{itemize} 
A. Blanchet and P. Lauren\c cot also proved in \cite{BL} the existence of self-similar compactly supported blowing-up solutions for $\overline{m}\in (M_c,M_c^+]$ where $M_c^+>M_c$.
See also \cite{Yao_Bertozzi} by Y. Yao and A.L. Bertozzi for recent formal and numerical results on self-similar and non self-similar blow-up for a generalization of system $(PKS^p)$ with kernel of power-law type.\\  

We can then observe similarities between both  problem $(PKS_{\frac{2}{N}})$ and $(PKS^{2-\frac{2}{N}})$ for $N\geq 3$. Here, we would like to thank P. Lauren\c cot for suggesting us that both problems should share the same stationary solutions, as we  can indeed verify, at least formally  : denoting $K_\Omega$ the Dirichlet kernel of the Laplacian in a bounded domain $\Omega$, it is easy to see 
that the steady states of $(PKS_q)$ and $(PKS^p)$ are respectively the solutions of 
$$ \frac{\rho^{1-q}}{1-q}-K_\Omega*\rho=C \mbox{ \qquad and \qquad } p\,\frac{\mu^{p-1}}{p-1} -K_\Omega*\mu =C',$$ 
where $C$ and $C'$ are any real constants.
Hence, the map $\rho \mapsto \mu:= (2-q)^{\frac{1}{q}}\, \rho$ defines a correspondence between steady states of $(PKS_q)$ with constant $C$ and steady states of $(PKS^{2-q})$ with constant $C'= (2-q)^{\frac{1}{q}}\, C$.\\
Then the formula for $M_c$ given  in \cite{BCL} also gives a formula for the critical cell mass $\overline{M}$ in our case :
$$\overline{M}=\left[ \frac{N^2\sigma_N}{C_*(N-1)}\right]^{\frac{N}{2}}\mbox{\quad  or equivalently \quad } M=\frac{N{\sigma_N}^{\frac{N}{2}-1}}{[C_*(N-1)]^{\frac{N}{2}}},$$ 
where $C_*$ is the optimal constant in the following variant of the Hardy-Littlewood-Sobolev inequality :  for all  $h\in
L^1(\mathbb{R}^N)\cap L^{2-\frac{2}{N}}(\mathbb{R}^N)$,
$$\left|\iint_{\mathbb{R}^N\times \mathbb{R}^N} \frac{h(x)h(y)}{|x-y|^{N-2}}dx\,dy \right|\leq C_*\|h\|_{2-\frac{2}{N}}^{2-\frac{2}{N}}\|h\|_1^{2/N} . $$
We would like to make the heuristic remark that if we roughly put $N=2$ in 
the above formula, then we recover the 
well-known critical mass $\overline{M}=8\pi$ since in this case $C_*=1$ and $\sigma_2=2\pi$.\\

It is interesting that, in spite of similarities in the qualitative behaviour, the two problems $(PKS_{\frac{2}{N}})$ and $(PKS^{2-\frac{2}{N}})$ 
seem to require different techniques. We would like to point out that our results are restricted to the radial setting, but on the other hand, they give a fairly more precise asymptotic description.\\

Throughout the rest of the paper, we assume : $$\boxed{m\geq 0,\, N\geq 3 \mbox{\, and \,} q=\frac{2}{N}.}$$

\section{The set of stationary solutions of problem $(PDE_m)$}

We begin by studying the steady states of problem $(PDE_m)$, i.e. the solutions of the following problem :
\begin{eqnarray}
x^{2-q}u''+u{u'}^q&=&0 \qquad x\in(0,1]\\
u(0)&=&0\\
u(1)&=&m
\end{eqnarray}
As is customary in evolution problems, this is essential in order to understand the large-time asymptotics of solutions of problem $(PDE_m)$.\\

Note that, even if we will use the results of this section only for $q=\frac{2}{N}$, they are all valid for any $q\in(0,1)$.

\subsection{Existence of a steady state depending on $m\geq 0$}

In view of proving Theorem \ref{thm_stat}, we first need to study the following Cauchy problem. 

\begin{Def} 
For $a\geq 0$, we define the problem $(P_a)$ by :
\begin{eqnarray}
\label{ode}
&& x^{2-q}u''+u{u'}^q=0 \\
\label{ci1}
&& u(0)=0  \\
\label{ci2}
&& u'(0)=a 
\end{eqnarray}
\end{Def}

\begin{Def} Let $R>0$.\\
We say that $u$ is a solution of problem $(P_a)$ on $[0;R[$ 
if :
\begin{itemize}
\item $u\in C^1([0;R[)\cap C^2(]0;R[)$.
\item $u$ is nondecreasing.
\item $u$ satisfies \eqref{ode} on $]0;R[$ and also conditions \eqref{ci1} and \eqref{ci2}.
\end{itemize}
This definition can obviously be adapted for the case of a closed interval $[0;R]$ or for $R=+\infty$.
\end{Def}

We first state some a priori properties of the solutions of $(P_a)$.

\begin{lem}
\label{lem1}
Let $a\geq 0$ and $R>0$.  \\
Let $u$ be a solution of $(P_a)$ on $[0;R[$. Then :
\begin{itemize}
\item [i)] If there exists $x_0\in[0;R[$ such that $u'(x_0)=0$ then $u(x)=u(x_0)$ 
for all $x\in [x_0;R[$.
\item [ii)] For all $x\in[0;R[$, $0\leq u(x) \leq 2$ and $0\leq u'(x) \leq a$.
\end{itemize}
 
\end{lem}

\textit{Proof :} 
i) We first note that since $u'\geq 0$, then, thanks to $(\ref{ode})$, $u$ is concave on $]0;R[$. If $x_0>0$ then for all $x\in[x_0;R[$, $0\leq u'(x)\leq u'(x_0)=0$. If $x_0=0$, we have to use in addition the continuity of $u'$ at $x=0$.\\
ii) If $a=0$ then by $i)$, $u=0$, hence the result.\\
Now, let us suppose $a>0$ and fix $x\in[0;R[$. 
$u'$ is nonincreasing on $]0;R]$ and continuous at $x=0$, so $0\leq u'(x) \leq a$. 
Hence, we just have to prove that $u(x)\leq 2$. \\
Let us denote $x_0=\sup\{x\in[0;R[, u'(x)>0 \}>0$. Let $x\in[0;x_0[$. \\
First note that by $i)$, $u'(x)>0$.
Since $u$ is concave, then $u'(x)\leq \dfrac{u(x)}{x}$ so 
$-\dfrac{u''(x)}{u'(x)^{q+1}}\geq \dfrac{1 }{x^{1-q}}$. Hence, by integration, we get
 $$\dfrac{1}{u'(x)^q}\geq \dfrac{1}{u'(x)^q}-\dfrac{1}{u'(0)^q}\geq x^q $$ 
and finally
 $xu'(x)\leq 1$. Since $u''(x)$ in nonpositive and  $$uu'(x)=-xu''(x)(xu'(x))^{1-q},$$ then 
$uu'(x)\leq -xu''(x)=(u-xu')'(x)$. Finally, by integration,  
$$\dfrac{u^2(x)}{2}\leq u(x)-xu'(x)\leq u(x)$$
so $u(x)\leq2$ for all $x\in[0;x_0[$.\\
If $x_0=R$, all is done. 
Else, $u'(x_0)=0$ and $u(x_0)\leq 2$ by continuity, so by $i)$, $u(x)=u(x_0)\leq 2$  
for all $x\in[x_0;R[$. \hfill $\square$\\

We will prove that solutions of problem $(P_a)$ exist on $[0,\infty)$.
We begin by showing the local existence.

\begin{lem}
\label{lem_existence_locale_sol_stationnaire} 
 Let $a\geq 0$ and $\tau>0$. \\
If $\tau$ is small enough, there exists a unique solution of $(P_a)$ on $[0;\tau]$.
\end{lem}

\textit{Proof :} If $a=0$, then from the previous lemma $i)$, it is clear that $0$ is the unique solution of the problem $(P_0)$  on $[0,\infty)$.\\
If $a>0$, 
let us define $$E_a=\{u\in C^1([0,\tau]),\;  u(0)=0,\; u'(0)=a,\;\|u'-a \|_\infty\leq
 \dfrac{a}{2} \}.$$ 
$E_a$ equipped with the metric induced by the norm $\|u\|_{E_a}=\|u'\|_\infty$ is 
a complete metric space. Any $u\in E_a$ is nondecreasing on $[0,\tau]$ since $u'\geq \frac{a}{2}$.\\
It is clear that the following function $F$ is well defined : $$F:E_a\rightarrow C^1([0;\tau])$$  $$F(u)(x)=ax-\int_0^x\int_0^y \dfrac{u(s)}{s}\dfrac{u'(s)^q}{s^{1-q}}ds
 dy.$$ 
Since for all $u\in E_a$, $\|u' \|_\infty \leq \dfrac{3}{2}a$, we easily get that  $$\|F(u)'-a\|_\infty  \leq (\frac{3}{2}a)^{q+1} \frac{\tau^q}{q}\leq \frac{a}{2}$$ if $\tau$ is chosen small enough. 
Hence, $F$ sends $E_a$ into $E_a$. \\
By the mean value theorem, since for all $u\in E_a$, $\frac{a}{2} \leq u'  \leq 
\frac{3}{2}a$,
 
$$\| F(u)-F(v)\|_{E_a}
\leq \left(\frac{(\frac{3}{2}a) ^q}{q}+\frac{\frac{3}{2}a}{(\frac{a}{2})^{1-q}}\right)\tau^q \|v-u\| _{E_a} 
$$ 
Hence, if $\tau$ is small enough, $F$ is a contraction so there exists a fixed point of $F$. Since 
$F(u)\in C^2((0,\tau])$ when $u\in C^1([0;\tau])$ then $u$ is a solution of 
$(P_a)$.\\
Finally, it is easy to check that a solution is necessarily a fixed point of $F$, which proves the uniqueness. \hfill $\square$

\begin{rmq}
Let $u$ be the local solution of $(P_a)$ on $[0,\tau]$. We would like to stress that $u$ is not $C^2$ up to $x=0$. Indeed, one can see that
$$ u'(x)=a-\frac{a^{q+1}}{q}x^q+o(x^q).$$
Moreover, one can obtain an expansion of $u$ at any order in powers of $x^q$. We proved in \cite{Montaru} that solutions of problem $(PDE_m)$ share these properties with the stationary solutions.
\end{rmq}

\begin{thm} 
\label{thm_existence_globale_sol_stationnaires}
Let $a\geq 0$.\\
\noindent There exists a unique maximal solution of $(P_a)$. \\Moreover, it is globally defined  on $[0,\infty)$.
\end{thm}

\textit{Proof :} For the sake of completeness, we prefer to give a (standard) proof. \\
\underline{Existence :} Leaving aside the obvious case $a=0$,
let $a>0$. \\
By Lemma \ref{lem_existence_locale_sol_stationnaire}, for a given $\tau$ small enough, we have a unique classical 
solution $u_\tau$ of $(P_a)$ on $[0;\tau]$. 
Setting $W=(u,u')$, we can now consider the following ordinary differential equation on the interval $[\dfrac{\tau}{2},+\infty[$ :
\begin{eqnarray}
\label{probleme_auxiliaire_edo}
 W'&=&F(x,W)\\
 W\left(\frac{\tau}{2}\right)&=&\left(u_\tau\left(\frac{\tau}{2}\right),u_\tau'\left(\frac{\tau}{2}\right)\right) \nonumber 
\end{eqnarray}
where $F(x,u,v)=\left(v,
\dfrac{-uv^q}{x^{2-q}}\right).$
Let us denote $\Omega=\mathbb{R} \times ]0;+\infty[$. Since $F$ is locally 
Lipschitz continuous with respect to $W$ in $\Omega$, by classical Cauchy-Lipschitz theory,  there exists a maximal solution  $u\in C^2([\frac{\tau}{2};X^*))$ of problem (\ref{probleme_auxiliaire_edo})
  such that $$(u(x),u'(x))\in \Omega \mbox{ for all } x\in
 [\frac{\tau}{2};X^*).$$ By local uniqueness in the classical Cauchy-Lipschitz theory, $u_\tau=u$ around $x=\frac{\tau}{2}$, so we can glue $u_\tau$ and $u$ and   get a solution of problem $(P_a)$ on $[0;X^*)$.\\
 If $X^*=+\infty$, all is done. So we suppose that $X^*<\infty$.\\
 Since $u$ is nondecreasing and bounded above by Lemma \ref{lem1} ii), then 
 $$l=\underset{x \rightarrow X^*}{\lim}u(x) \mbox{ exists.}$$ 
 Hence, we can extend continuously the function $u$ by 
setting $u(x)=l$ for $x\geq X^*$.
But $(u(x),u'(x))$ must leave any compact of $\Omega$ as $x$ goes to $X^*$ so, by 
 Lemma $\ref{lem1})$ ii), the only possibility is that $\underset{x \rightarrow X^*}{\lim}u'(x)=0$.
So $u\in C^1([0;+\infty))$. \\
And now, thanks to 
$(\ref{ode})$, $\underset{x \rightarrow X^*}{\lim}u''(x)=0$, so  $u\in C^2((0;+\infty)$. Moreover, $u$ clearly satisfies $(\ref{ode})$ on $(0,+\infty)$ and is then  a global solution of $(P_a)$.\\
\underline{Uniqueness :} Let v another global solution of $(P_a)$. By the result of uniqueness around $x=0$ and the uniqueness due to classical Cauchy-Lipschitz theory in $\Omega$, $u$ and $v$ coincide for all $x\in [0,X^*)$.
If $X^*=\infty$, all is done.  Now, assume that $X^*<+\infty$.   $v(X^*)=u(X^*)$  by continuity. As $v'(X^*)=u'(X^*)=0$, then by Lemma $\ref{lem1})$ i), $v(x)=v(X^*)=u(X^*)$ for all $x\geq X^*$. Hence, $v=u$. \hfill $\square$

\begin{nota}
Let $a \geq 0$. \\ 
We denote $U_a$ the unique solution of $(P_a)$ on $ [0,+\infty)$.
\end{nota}

\begin{lem}
Let $a\geq 0$. \\
There exists $x_0\geq 0$ such that $U_a(x)=U_a(x_0)$ for all $x\geq x_0$. 
\end{lem}

\textit{Proof :} Suppose the contrary. By Lemma $\ref{lem1})i)$, it implies that $U_a'(x)>0$ for all $x\geq 0$. Since $U_a$ is nondecreasing and has 2 as an upper bound, there exists $l\leq 2$ such that $U_a(x)$ tends to $l$ as $x$ goes to infinity. \\
As $U_a'$ is nonnegative and nonincreasing, then $U_a'(x)$ has a nonnegative limit as $x$ goes to $+\infty$, but this limit has to be $0$ since $U_a$ is bounded from above. \\
Moreover, for all $x>0$,
$$\dfrac{d}{dx}U_a'(x)^{1-q}=U_a(x)\dfrac{d}{dx}\dfrac{1}{x^{1-q}}. $$
Let $x_0>0$ and $x\geq x_0$. By monotonicity, $U_a(x_0)\leq U_a(x)\leq l$ so
by integration on $[x_0;x]$,
$$ l \left( \frac{1}{x^{1-q}} -\frac{1}{x_0 ^{1-q}}\right)  \leq U_a'(x)^{1-q}-U_a'(x_0)^{1-q}\leq U_a(x_0) \left( \frac{1}{x^{1-q}} -\frac{1}{x_0 ^{1-q}}\right) . $$ 

Finally, we let $x$ go to $+\infty$ and then $x_0$ go to $+\infty$ to obtain that 
$$ U_a'(x_0) \underset{x_0 \rightarrow +\infty}{\sim} \dfrac{l^{\frac{1}{1-q}}}{x_0}. $$
Then, $U_a(x)$ goes to infinity as $x$ goes to $\infty$, which is a contradiction. \hfill $\square$

\begin{propo}
$U_a(x)=U_1(ax)$ for all $x\geq 0$ and all $a \geq 0$.
\end{propo}

\textit{Proof :} Let $V(x)=U_1(ax)$. Clearly, $V\in C^1([0,+\infty))\cap C^2((0;+\infty))$, $V(0)=0$, $V'(0)=a$ and for $x>0$, $$x^{2-q}V''(x)=(ax)^{2-q}U_1''(ax)a^q=U_1(ax)U_1'(ax)^q a^q=V(x)V'(x)^q.$$
The result then follows from the uniqueness of the solution of problem $(P_a)$. \hfill $\square$

\begin{rmq}
Behind this proof is the fact that $L_a$ and $D$ commute, where for any 
$u\in C^1([0,\infty))$, $a\geq 0$ and $x\geq 0$, $ L_a(u)(x)=u(ax)$ and $Du=x\,u_x$.
\end{rmq}

This proposition drives us to the natural following definition.
\begin{Def}
The number $M=\underset{x\geq 0}{\max} \;U_1(x)$ will be called the critical mass.\\
Note that $M$ also is the maximal value of each $U_a$, for all $a>0$. 
\end{Def}

We also will use the following notation.
\begin{nota}
We denote by $A$ the first $x\geq 0$ such that $U_1(x)=M$.
\end{nota}

\begin{propo} \quad
\begin{itemize}
\item[i)] If $m \in [0;M)$, there exists a unique $a\in[0;A)$ such that $U_a(1)=m$.
\item [ii)] $U_a(1)=M$ if and only if $a\geq A$.
\end{itemize}
\end{propo}

\textit{Proof :} 
i) $U_a(1)=U_1(a)$ and $U_1$ is a bijection from $[0;A)$ to $[0;M)$.\\
ii) $U_a(1)=M$ if and only if $U_1(a)=M$, which is equivalent to $a\geq A$. \hfill $\square$ \\

From this follows Theorem \ref{thm_stat} announced in the introduction.

\subsection{Order and topological properties of the set of stationary solutions}

Now, we shall describe two simple but very important properties of the set of stationary solutions : 
it has a total order and its topology behaves well with respect to that order.

\begin{nota}
Let $m\geq 0$. 
\begin{itemize}
\item $Y_m^1=\{u\in C^1([0;1]),\: u \: nondecreasing,\, u(0)=0,\,u(1)=m\}$ is the complete metric space equipped with the distance induced
by the $C^1$ norm. 
\item $V_m(u,\epsilon)\subset Y^1_m$ is the open ball of center $u\in Y^1_m$ and radius $\epsilon > 0$.
\item We say that two functions $u,v\in Y^1_m$ satisfy $u\prec v$ if  $u\leq v$ and $u \neq v$.\\
If $u\in Y^1_m$ and $V\subset Y^1_m $, we say that $u\prec V$ if for all $v\in V$, $u\prec v$. \\
Similarly, we define $V\prec u$.
\end{itemize}
\end{nota}

\begin{propo} \quad
\label{prop_stat}
\begin{itemize} 
 \item [i)] Suppose $0\leq a< b$. Then $U_a\prec U_b$.\\
 The set of stationary solutions is a totally ordered set.
\item [ii)] Suppose $A\leq a<b$.\\
 $\alpha)$ There exists $\epsilon>0$ such that 
$V_M(U_a,\epsilon) \prec U_b$.\\ 
$\beta$) There exists $\epsilon>0$ such that 
$\{u\in Y_M^1, \; u \leq U_a\}\cap V_M(U_b,\epsilon)= \emptyset$.

\end{itemize}
\end{propo}

\textit{Proof :} i) Let $x\geq 0$. $U_a(x)=U_1(ax)\leq U_1(bx)=U_b(x)$, so $U_a\leq U_b$. \\
Since, $U_a'(0)<  U_b'(0)$, then $U_a\prec U_b$.\\
ii) $\alpha)$ There exists $\gamma>0$ such that for all $x\in [0;\gamma]$, $ \frac{a+b}{2}x \leq U_b(x)$.\\
Let us set $\epsilon_1=\frac{b-a}{2}$, $\epsilon_2=\frac{1}{2}\underset{x\in [\gamma;\frac{A}{b}]}
{\min} \left[ U_b(x)-U_a(x) \right]$ and $\epsilon=\min(\epsilon_1,\epsilon_2)$. \\
Note that $\epsilon>0$ since for all $x\in [\gamma;\frac{A}{b}]$, we have 
$ax<bx\leq A$, hence $$U_a(x)=U_1(ax)<U_1(bx)=U_b(x)$$ because $U_1$ is increasing on $[0;A]$.\\
Let $u\in V_M(U_a,\epsilon)$. 
Since $u\in V_M(U_a,\epsilon_1)$, then $$u'\leq U_a'+\frac{b-a}{2}\leq a+\frac{b-a}{2}=\frac{a+b}{2}.$$
Hence, for all $x\in [0;\gamma]$, $$ u(x)\leq \frac{a+b}{2}x \leq U_b(x).$$
Since $u\in V_M(U_a,\epsilon_2)$, it is clear that for all $x\in [\gamma;\frac{A}{b}]$, $u(x) < U_b(x).$
Moreover, if $x\geq \frac{A}{b}$ then $u(x)\leq M=U_b(x)$. Hence, $$u\prec U_b.$$
$\beta)$ Let $u \in Y_M^1$ such that $u \leq U_a$. Then $u'(0)\leq a$. 
Let us set $\epsilon=\frac{b-a}{2}>0$.
Since $\|U_b-u\|_{C^1}\geq U_b'(0)-u'(0)\geq b-a$ then $u \notin V_M(U_b,\epsilon)$. \hfill $\square$

\section{Summary of local in time results}
\label{prelim}

In this section, we give some useful results on problem $(PDE_m)$ and two other auxiliary 
parabolic problems.
For proofs, see \cite{Montaru}.

\subsection{Wellposedness of problem $(PDE_m)$}
\label{problem_PDE}
Before stating 
the existence and uniqueness of classical solutions for problem $(PDE_m)$, 
we need to fix some definitions.

\begin{Def} 
$$Y_m=\{u\in C([0;1]),\:u \mbox{ nondecreasing },\; u'(0) \mbox{ exists, } u(0)=0,\;u(1)=m\}.$$
\end{Def}

\begin{Def} Let $T>0$. \\
We say that $u$ is a classical solution of $(PDE_m)$ (see (\ref{equ_u_1})) on $[0,T)$ with initial condition $u_0\in Y_m$ if :
\begin{itemize}
\item $u\in C([0,T)\times [0,1])\bigcap
C^1((0,T)\times [0,1])\bigcap C^{1,2}((0,T)\times (0,1]).$
\item $u(0)=u_0.$ 
\item $u(t)\in Y_m$ for $t\in [0,T).$
\item  $u$ satisfies $(\ref{equ_u_1})$ on $(0,T)\times (0,1]$.
\end{itemize}
\end{Def}

\begin{Def} For any real function defined on $(0,1]$, we set :
$$\mathcal{N}[u]=\underset{x\in(0,1]}{\sup}\frac{u(x)}{x}.$$
\end{Def}

\begin{thm}Let $q\in(0,1)$, $m\geq 0$ and $u_0\in Y_m $.
\label{thm_existence_u}
\begin{itemize}
 \item [i)]
There exists $T_{max}=T_{max}(u_0)>0$ and a unique maximal classical solution $u$ of problem $(PDE_m)$ with initial condition $u_0$. \\
Moreover, $u$ satisfies the following condition :
\begin{equation}
\label{maj_u_C1}
 \underset{t\in(0,T]}{\sup} \sqrt{t}\; \|u(t)\|_{C^1([0,1])}<\infty \mbox{ for any } T\in(0,T_{max}).
\end{equation}
\item[ii)] Blow-up alternative : 
$T_{max}=+\infty$ \quad or \quad  $\underset{t\rightarrow T_{max}}
{\lim}\mathcal{N}[u(t)]=+\infty .$
\item[iii)] $u_x(t,0)>0$ for all $t\in (0,T_{max})$.
\end{itemize}
\end{thm}

Moreover, a classical comparison principle is available for problem $(PDE_m)$. 
 
\begin{lem}
\label{PC_PDE}
Let $T>0$. Assume that :
\begin{itemize}
\item $u_1,u_2\in C([0,T]\times [0,1])\bigcap
C^1((0,T]\times [0,1])\bigcap C^{1,2}((0,T]\times (0,1)).$ 
\item For all $t\in (0,T] $, $u_1(t)$ and $u_2(t)$ are nondecreasing. 
\item There exists $i_0\in\{1,2\}$ and some $\gamma<\frac{1}{q}$ such that 
$$ \underset{t\in(0,T]}{\sup} t^\gamma\; \|u_{i_0}(t)\|_{C^1([0,1])}<\infty .$$
\end{itemize}
Suppose moreover that :
\begin{eqnarray} 
 & (u_1)_t\leq x^{2-\frac{2}{N}}(u_1)_{xx}+u_1 (u_1)_x^q  & \mbox{ for all }(t,x)\in (0,T]\times (0,1)\\
 & (u_2)_t\geq x^{2-\frac{2}{N}}(u_2)_{xx}+u_2 (u_2)_x^q  & \mbox{ for all }(t,x)\in (0,T]\times (0,1)\\
& u_1(0,x)\leq u_2(0,x) &\mbox{ for all }x\in [0,1]\\
& u_1(t,0)\leq u_2(t,0) &\mbox{ for }t\geq 0 \\
& u_1(t,1)\leq u_2(t,1) &\mbox{ for }t\geq 0
\end{eqnarray}
Then $u_1\leq u_2$ on $[0,T]\times [0,1]$.\\
\end{lem}

\subsection{Two auxiliary parabolic problems}
\subsubsection{Problem $(tPDE_m)$}
\label{problem_tPDE}
We now introduce an auxiliary transformed problem $(tPDE_m)$ which will be helpful in order to get some estimates implying the compactness of the trajectories in the subcritical and critical case $m\leq M$.
This transformation was also important in \cite{Montaru} in order to establish the blow-up alternative (see Theorem \ref{thm_existence_w} ii) below), a property which will be used in the global existence part of the proof of Theorem \ref{thm_convergence}.\\

Denoting $B$ the open unit ball in 
$\mathbb{R}^{N+2} $ and $Z_m=\{w\in C(\overline{B}), \; w| _ { \partial B } = m\}$, we define 
$\begin{array}{lll}
\theta_0 &:& Y_m\longrightarrow Z_m\\
&& u_0\longmapsto w_0
\end{array}$
 where 
 $\left\{
\begin{array}{lll}
w_0(y)&=  \frac{u_0(|y|^N)}{|y|^N}&\mbox{ if } y\in 
\overline{B}\backslash \{0\} \\
&= u_0'(0) &\mbox{ if } y=0\\
\end{array}
\right. .$

We now make the following change of unknown 
\begin{eqnarray}
w(t,y)&=\frac{u(N^2t,|y|^N)}{|y|^N} &\mbox{ if } y\in 
\overline{B}\backslash \{0\}
\end{eqnarray}

in problem $(PDE_m)$ with initial condition $u_0\in Y_m$
and obtain the following problem $(tPDE_m)$ with initial condition $w_0=\theta_0(u_0)$.

\begin{Def}Let $m\geq0$ and $T>0$. \\
Let $u_0\in Y_m$ and $w_0=\theta_0(u_0)$.\\
We define problem $(tPDE_m)$ with initial condition $w_0$ by :
\begin{equation}
(tPDE_m) \left\lbrace \qquad 
\label{equ_w1}
\begin{array}{ccc}
 & w_t=\Delta w + N^2 w (w+\frac{y.\nabla 
w}{N})^q  &\mbox{ on }(0,T]\times \overline{B}\\
&w(0)=w_0&\\
&w+\frac{y.\nabla w}{N}\geq 0 & \mbox{ on }(0,T]\times \overline{B}\\
&w=m&\mbox{ on } [0,T]\times \partial B
\end{array}
\right.
\end{equation}

 A classical solution on $[0,T]$ for problem $(tPDE_m)$ with initial condition $w_0$  is a function 
$$ w\in C([0,T]\times \overline{B}) \bigcap C^{1,2} ((0,T]\times \overline{B})$$
such that all conditions of 
 $(\ref{equ_w1})$ are satisfied.\\
We define analogously a classical solution on $[0,T)$.
\end{Def}

Let us give now the corresponding wellposedness result.

\begin{thm} 
\label{thm_existence_w} Let $u_0\in Y_m$ and $w_0=\theta_0(u_0)$.
\begin{itemize}
 \item [i)] There exists $T^*=T^*(w_0)>0$ and a unique maximal 
classical radially symmetric solution $w$ of problem $(tPDE_m)$ with initial condition $w_0$.\\
Moreover, $w$ satisfies the following condition :
\begin{equation}
\label{majoration_w_C1}
 \underset{t\in(0,T]}{\sup} \sqrt{t}\; \|w(t)\|_{C^1(\overline{B})}<\infty \mbox{ for any } T\in(0,T^*).
\end{equation}
\item [ii)] Blow-up alternative : $T^*=+\infty \quad$ or $\quad\underset{t\rightarrow T^* }{\lim}
\|w(t)\|_{\infty, \overline{B}}=+\infty.$

\end{itemize}
\end{thm}

\textbf{Connection between problems $(PDE_m)$ and $(tPDE_m)$} :\\
Let $w_0=\theta_0(u_0)$ with $u_0\in Y_m$. Then,
 $$T_{max}(u_0)=N^2 T^*(w_0).$$
Moreover, for all $(t,x)\in [0,T_{max})\times [0,1]$, 
\begin{equation}
\label{def_u}
u(t,x)=x\; \tilde{w}(\frac{t}{N^2},x^{\frac{1}{N}}),
\end{equation}
 where we write $f=\tilde{f}(|\cdot|)$ for any radial function $f$ on $B$.

\subsubsection{Problem $(PDE_m^\epsilon)$}
\label{problem_PDE_eps}
We define an approximate problem $(PDE_m^\epsilon)$ of problem $(PDE_m)$.
It will be useful since we can easily find a strict Lyapunov functional $\mathcal{F}_\epsilon$ for $(PDE_m^\epsilon)$ and then prove by this way that $\mathcal{F}=\underset{\epsilon \rightarrow 0}{\lim}\,\mathcal{F}_\epsilon$ is a strict Lyapunov functional for the dynamical system induced by problem $(PDE_m)$ in the subcritical and critical case $m\leq M$.

\begin{Def}
Let $\epsilon>0$. We set :
$$f_\epsilon(x)=(x+\epsilon)^q-\epsilon^q \mbox{ \quad for }x\geq 0 $$
\end{Def}
Observe in particular that $0\leq f_\epsilon(x)\leq x^q$ for all $x\in[0,+\infty)$.

\begin{Def} Let $\epsilon>0$, $m\geq 0$ and $T>0$.\\
 We define problem $(PDE_m^\epsilon)$ with initial condition $u_0\in Y_m$ by :
 \begin{equation}
 (PDE_m^\epsilon) \left\lbrace \qquad
 \begin{array}{ccc}
\label{equ_u_eps1}
 & u_t=x^{2-\frac{2}{N}}u_{xx}+u f_\epsilon(u_x) &\mbox{ on } (0,T]\times (0,1] \\
& u(0)=u_0 & \\
&u(t,0)=0 &\mbox{ for all } t\in [0,T]\\
&u(t,1)=m &\mbox{ for all } t\in [0,T]
\end{array} 
\right.
 \end{equation}

 A classical solution on $[0,T]$ of problem $(PDE_m^\epsilon)$ with initial condition \\
 $u_0\in Y_m$ is a function 
$$u^\epsilon\in C([0,T]\times [0,1])\bigcap
C^1((0,T]\times [0,1])\bigcap C^{1,2}((0,T]\times (0,1])$$
such that all conditions of 
 $(\ref{equ_u_eps1})$ are satisfied.\\
A classical solution of problem $(PDE_m^\epsilon)$ on $[0,T)$ is defined analogously.
\end{Def}

\begin{thm}Let $m\geq 0$, $\epsilon>0$, $K>0$ and $u_0\in Y_m $ with $\mathcal{N}[u_0]\leq K$.
\label{thm_existence_u_epsilon}
\begin{itemize}
 \item [i)]
There exists $T_{max}^\epsilon=T_{max}^\epsilon(u_0)>0$ and a unique maximal 
classical solution $u^\epsilon$ on $[0,T_{max}^\epsilon)$ of  problem
$(PDE_m^\epsilon)$ with initial condition $u_0$.\\
Moreover, $u^\epsilon$ satisfies the following condition :
\begin{equation}
\label{maj_u_eps}
 \underset{t\in(0,T]}{\sup} \sqrt{t}\; \|u^\epsilon(t)\|_{C^1([0,1])}<\infty \mbox{ for all } T\in(0,T_{max}^\epsilon).
\end{equation}
\item[ii)] Blow-up alternative : 
$T^\epsilon_{max}=\infty$ \quad or \quad  $\underset{t\rightarrow T_{max}^\epsilon}
{\lim}\mathcal{N}[u^\epsilon(t)]=+\infty .$
\item[iii)] $(u^\epsilon)_x>0$ on $t\in (0,T_{max}^\epsilon)\times [0,1]$.
\item [iv)] $u^\epsilon\in C^2\left((0,T_{max}^\epsilon)\times (0,1]\right)$. 
\textit{(not optimal)}
\end{itemize}
\end{thm}

\textbf{Connection with problem $(PDE_m)$ :} \\
Let us fix an initial condition $u_0\in Y_m$.\\ 
The next lemma shows the convergence of maximal classical solutions $u^\epsilon$ of $(PDE_m^\epsilon)$ to the maximal classical solution of $(PDE_m)$ in various spaces. \\
These results are essential in our proof that $\mathcal{F}=\underset{\epsilon \rightarrow 0}{\lim}\,\mathcal{F}_\epsilon$ is a strict Lyapunov functional in the case $m\leq M$.

\begin{lem}
\label{lem_convergence_u_epsilon}
Let $u_0\in Y_m$.
\begin{itemize}
\item[i)] $T_{max}(u_0)\leq T_{max}^\epsilon(u_0)$ for any $\epsilon>0$.
\item[ii)]  Let $[t_0,T]\subset (0,T_{max}(u_0))$.
\begin{itemize}
\item [$\alpha)$]$u^\epsilon \underset{\epsilon \rightarrow 0}{\longrightarrow} u$ in $C^{1,2}([t_0,T]\times (0,1]).$\\
Moreover, there exists $K>0$ independent of $\epsilon$ such that\\
for all $(t,x)\in[t_0,T]\times (0,1]$, $|u^\epsilon_{xx}|\leq \frac{K}{x^{1-q}}$.
\item[$\beta )$] $(u^\epsilon)_x \underset{\epsilon \rightarrow 0}{\longrightarrow} u_x$ in $C([t_0,T]\times [0,1])$.
\item[$\gamma )$] $(u^\epsilon)_t \underset{\epsilon \rightarrow 0}{\longrightarrow} u_t$ in $C([t_0,T]\times [0,1])$.
\end{itemize}
\end{itemize} 
\end{lem}

\section{Convergence to a stationary state in critical and subcritical case $m\leq M$}

All this section only concerns the case $m\leq M$.\\

We shall prove that problem $(PDE_m)$ defines a continuous dynamical system on $Y_m^1$  which admits a strict Lyapunov functional. We shall be able to prove that classical solutions
of $(PDE_m)$ converge to a stationary state as times goes to infinity, even in the case $m=M$ where there is a continuum of steady states. 

\subsection{Estimates}

\begin{lem}
\label{traj_globales}
 Let $m\leq M$ and $u_0\in Y_m$. 
 Then, $$T_{max}(u_0)=\infty.$$ 
 Moreover, for each $K>0$,  for any $u_0\in Y_m$ with  $\mathcal{N}[u_0]\leq K$, we have  
$$\underset{t\in [0,\infty)}{\sup}\mathcal{N}[u(t)]\leq C_K$$
where $C_K=\frac{A}{M}\,\max(K,M)$.
\end{lem}

\textit{Proof : } Let $T_{max}=T_{max}(u_0)$.\\
From Theorem $\ref{thm_existence_u}$, it is sufficient to prove that 
$$\underset{t\in[0,T_{max})}{\sup}\mathcal{N}[u(t)]<\infty .$$
This fact  easily follows from a comparison with a supersolution of problem $(PDE_m)$.
The main idea is that since $m\leq M$, if $a$ is large enough then $u_0\leq U_a$ and $U_a$ is then 
a supersolution so $0\leq u(t)\leq U_a$ for all $t\in [ 0,T_{max})$. \\
More precisely, since $u_0$ is differentiable at $x=0$, $x\mapsto \frac{u_0(x)}{x}$ can be extended continuously to $[0;1]$, so there exists 
$K\geq M$ such that 
$\mathcal{N}[u_0] \leq K$.\\
Let us set $a= \frac{K}{M}A$. \\
For $x\in [0;\frac{M}{K}]$, $u_0(x)\leq Kx\leq U_a(x)$ since by concavity, $U_a$ is above its chord 
between $x=0$ and $x=\frac{M}{K}$. On $[\frac{M}{K},1]$, $u_0\leq m\leq  M=U_a$. \\
Hence, $u_0\leq 
U_a$ on $[0,1]$. Finally, since $M\geq m$, $U_a$ is a supersolution of the PDE, so $u(t)\leq U_a$ for all $t\in [ 0,T_{max})$.
By concavity of $U_a$, we see that $ \mathcal{N}[U_a]=a$, so
$\underset{t\in[0,T_{max})}{\sup}\mathcal{N}[u(t)]\leq a<\infty$. 
Then $T_{max}=\infty$.\\
We notice that the choice of $a$ depends only on $K$, whence the second part of the lemma.\hfill $\square$\\

Before going further, we would like to recall some notation and properties of the heat semigroup.
For reference, see for instance the book \cite{Lu} of A. Lunardi.

\begin{nota} \quad 
\begin{itemize}
\item $B$ denotes the open unit ball in $\mathbb{R}^{N+2}$.
\item $X_0=\{W \in C(\overline{B}), \; 
W | _{\partial B}=0 \}$.
\item $(S(t))_{t\geq 0}$ denotes the heat semigroup  on $X_0$. It is the restriction on $X_0$ of the Dirichlet heat semigroup on $L^2(B)$.
\item $(X_\theta)_{\theta \in [0,1]}$ denotes the scale of interpolation spaces for 
$(S(t))_{t\geq 0}$, where $X_0=L^2(B)$, $X_1=D(-\Delta)$ and $X_\alpha \hookrightarrow X_\beta$ with dense continuous injection for any $\alpha>\beta$, $(\alpha,\beta)\in[0,1]^2$.
\end{itemize}

\end{nota}

\begin{prop}\quad 
\begin{itemize}
\item $X_{\frac{1}{2}}=\{ W \in C^1(\overline{B}), \;W | _{\partial B}=0 \} . $
\item Let $\gamma_0 \in (0;\frac{1}{2}]$. For any $\gamma\in[0,2\gamma_0)$,  $$X_{\frac{1}{2}+\gamma_0} \subset C^{1,\gamma}
(\overline{B})$$  with continuous embedding.
\item There exists $C_D\geq 1$ such that for any $\theta \in [0;1]$,
 $W \in C(\overline{B})$ and $t>0$,
$$\|S(t)W \|_{X_\theta}\leq \frac{C_D}{t^{\theta}}\|W\|_{\infty}. $$
\end{itemize}

\end{prop}

We just want to introduce a specific notation we are going to use.
\begin{nota} 
Let $(a,b)\in(0,1)^2$. We denote $I(a,b)=\int_0^1 \frac{ds}{(1-s)^as^b}$.  \\ For all $ t\geq 0,\;
\int_0^t \frac{ds}{(t-s)^as^b}=t ^{1-a-b}I(a,b) .$

\end{nota}

We will now give an estimate from which follows a compactness result. 

\begin{lem}
\label{lem_compacite}
 Let $m\leq M $, $\gamma\in [0;1)$, $t_0>0$ and $K>0$.\\
There exists $D_K>0$ such that for any $u_0\in Y_m $ with $\mathcal{N}[u_0]\leq K$, 
then $$\underset {t \geq t_0}{\sup \;} \| 
 u(t) \| _{C^{1,\frac{\gamma}{N}}} \leq D_K.$$
 As a consequence, $\{u(t), \;\mathcal{N}[u_0]\leq K,\; t\geq t_0\}$ is relatively 
compact in $Y_m^1$.
\end{lem}

\textit{Proof :} Let $u_0\in Y_m $ such that $\mathcal{N}[u_0]\leq K$. 
Let $w_0=\theta_0(u_0)$.\\
\underline{First step :} from Lemma \ref{traj_globales}, there exists $C_K>0$ such that 
$$\underset{t\in [0,\infty)}{\sup}\mathcal{N}[u(t)]\leq C_K.$$
Since for $t\geq0$, $\|w(t)\|_{\infty,\overline{B}}=\mathcal{N}[u(\frac{t}{N^2})]$, we 
deduce that $w$ is global and that 
$$\underset{t\in [0,\infty)}{\sup}\|w(t)\|_{\infty,\overline{B}}=\underset{t\in [0,\infty)}
{\sup}\mathcal{N}[u(t)]\leq C_K.$$

 \underline{Second step :} Let $\tau=\frac{t_0}{N^2} $ and $t\in [0,\tau]$. \\
Denoting $W_0=w_0-m$, then
 \begin{equation}
 \label{w_mild}
 w(t)-m=S(t)W_0+\int_0^t S(t-s)N^2w \left(w+\frac{x.\nabla 
 w}{N}\right)^q ds,
 \end{equation}
 so
 $$ \|w(t)\|_{C^1}\leq m+ \frac{C_D}{\sqrt{t}}(C_K+m) +N^2\int_0^t \frac{C_D}{\sqrt{t-s}}C_K 
 \|w(s)\|_{C^1}^q ds .$$
 Setting $h(t)=\underset{s\in (0,t]}{\sup} \sqrt{s}\|w(s)\|_{C^1}$, we have $h(t)<\infty$ by (\ref{majoration_w_C1}) and 
 $$ \sqrt{t}\|w(t)\|_{C^1}\leq m \sqrt{\tau}+ C_D(C_K+m) +N^2C_KC_D\sqrt{t}\int_0^t \frac{1}
 {s^{\frac{q}{2}}\sqrt{t-s}} h(s)^q ds ,$$
 $$ \sqrt{t}\|w(t)\|_{C^1}\leq m \sqrt{\tau}+ C_D(m+C_K) +N^2C_KC_D I\left(\frac{1}{2},\frac{q}{2}\right) 
 t^{1-\frac{q}{2}} h(t)^q  .$$
 Let $T\in(0,\tau]$. Then,
 \begin{equation}
 \label{h(T)}
  h(T)\leq m \sqrt{\tau}+ C_D(m+C_K) +N^2C_KC_D I\left(\frac{1}{2},\frac{q}{2}\right) T^{1-\frac{q}{2}} h(T)^q .
   \end{equation}
Setting $A=m \sqrt{\tau}+ C_D(m+C_K)$ and $B=N^2C_KC_D I \left( \frac{1}{2},\frac{q}{2} \right)  2^q$, assume that 
there exists $T\in [0,\tau]$ such that $h(T)=2A$. Then,
$$A^{1-q}\leq B T^{1-\frac{q}{2}} \mbox{ which implies } T\geq \left(\frac{A^{1-q}}{B}\right)^
{\frac{1}{1-q}}.$$
Let us set $\tau'=\min\left(\tau,\frac{1}{2}\left(\frac{A^{1-q}}{B}\right)^{\frac{1}{1-q}} \right)$. \\
Since $h\geq 0$ is nondecreasing, $h_0=\underset{t\rightarrow 0^+}{\lim} h(t)$ exists and  $h_0\leq A$ by (\ref{h(T)}). So by continuity of $h$ on $(0,\tau']$, 
$ h(t)\leq 2A \mbox{ for all }t\in(0,\tau'] $, that is to say :
$$ \|w(t)\|_{C^1}\leq \frac{2A}{\sqrt{t}} \mbox{ for all }t\in(0,
\tau'],$$ 
where $A$ and $\tau'$ only depend on $K$. Then, setting $A_K=2A$, we have
$$\underset{t\in[0,\tau']}{\sup}\sqrt{t}\|w(t)\|_{C^1}\leq A_K.$$

\underline{Third step :}
Let $\gamma_0\in (\frac{\gamma}{2},\frac{1}{2})$ and $t\in[0,\tau']$.\\
Setting $W=w-m$ and $W_0=w_0-m$, then for $t\geq 0$, due to (\ref{w_mild}), we get 
\begin{align*} 
  \|W(t)\|_{X_{\frac{1}{2}+\gamma_0}}&\leq  \frac{C_D}{t^{\frac{1}{2}+\gamma_0}}(C_K+m) +N^2\int_0^t 
\frac{C_D}{(t-s)^{\frac{1}{2}+\gamma_0}}C_K 
 \frac{(A_K)^q}{s^{\frac{q}{2}}} ds .&
   \end{align*}
 Then we deduce that :
 \begin{align*} 
  t^{\frac{1}{2}+\gamma_0}\|W(t)\|_{X_{\frac{1}{2}+\gamma_0}}&\leq  C_D(C_K+m) +N^2C_K C_D (A_K)^q t^{\frac{1}{2}+\gamma_0}\int_0^t \frac{1}
 {(t-s)^{\frac{1}{2}+\gamma_0} s^{\frac{q}{2}}}  ds &\\
  &\leq C_D(m+C_K) +N^2C_K C_D (A_K)^q I(\frac{1}{2}+\gamma_0,\frac{q}{2}) 
 {\tau'}^{(1-\frac{q}{2})}.&  
  \end{align*}
Hence, since $X_{\frac{1}{2}+\gamma_0}\subset C^{1,\gamma}(\overline{B})$, we deduce that there exists $A_K'>0$ depending only on $K$ such that  
$\underset{t\in[0,\tau']}{\sup}t^{\frac{1}{2}+\gamma}\|w(t)\|_{C^{1,\gamma}(\overline{B})}\leq A_K'$. Then, 
$$\|w(\tau')\|_{C^{1,\gamma}(\overline{B})} \leq \frac{A_K'}{\tau'^{(\frac{1}{2}+\gamma)}}=:A_K''.$$
\underline{Last step :}
Let $t'\geq \frac{t_0}{N^2}$. Since $\tau'\leq\frac{t_0}{N^2}$, we can apply the same arguments by taking $w_0(t'-\tau)$ as initial data instead of $w_0$, so 
we obtain $$\mbox{for all } t'\geq \frac{t_0}{N^2},\,\|w(t')\|_{C^{1,\gamma}(\overline{B})} \leq A_K'' .$$
Finally, coming back to $u(t)$, thanks to formula (\ref{def_u}), we get an upper bound $D_K$ for $ \|  u(t) \| _{C^{1,\frac{\gamma}{N}}}$ valid for any $u_0\in Y_m$ such that $\mathcal{N}[u_0]\leq K$. \hfill $\square$

\subsection{A continuous dynamical system $(T(t))_{t\geq 0}$}

We recall the definition of a continuous dynamical system on $Y_m^1$.\\
For reference, see \cite[chap. 9, p.142]{Haraux}.
\begin{Def}
A continuous dynamical system on $Y_m^1$ is a one-parameter family of mappings $(T(t))_{t\geq 0}$ from $Y_m^1$ to $Y_m^1$   such that :
\begin{itemize}
\item [i)] $T(0)=Id$.
\item [ii)]$T(t+s)=T(t)T(s)$ for any $t,s\geq 0$.
\item [iii)]For any $t\geq 0$, $T(t)\in C(Y_m^1,Y_m^1)$.
\item [iv)]For any $u_0\in Y_m$, $t\mapsto T(t)u_0 \, \in C((0,\infty),Y_m^1)$.
\end{itemize}
\end{Def}

\begin{rmq}
Continuity at $t=0$ is sometimes included in the definition, but it is not required for our needs.
\end{rmq}

\begin{Def}
Let $u_0\in Y_m^1$ and $t\geq 0$.\\ 
We define $T(t)u_0=u(t)$ where we recall that $u$ is the classical solution of problem $(PDE_m)$ with initial condition $u_0$.
\end{Def}

\begin{propo}
$(T(t))_{t\geq 0}$ is a continuous dynamical system on $Y_m^1$.
\end{propo}

\textit{Proof : } Thanks to Lemma \ref{traj_globales}, we know that $T(t)$ 
is well defined for all $t\geq 0$ and by definition of a classical solution, $T(t)$ maps $Y_m^1$ into $Y_m^1$. \\
ii) is clear by uniqueness of the global classical solution.\\
iv) comes from the fact that $u\in C((0,\infty),C^1([0,1]))$.\\
iii) Let $t> 0$, $u_0\in Y_m^1$ and $(u_n)_{n\geq 1}\in Y_m^1$.\\
 Assume that $u_n \underset{n \rightarrow \infty}{\overset{C^1}{\longrightarrow}} u_0$.
Let us show that $u_n(t) \underset{n \rightarrow \infty}{\overset{C^1}{\longrightarrow}} u(t)$.\\
We proceed in two steps.\\
\underline{First step :}  We show that if $u_n \underset{n \rightarrow \infty}{\overset{C^0}{\longrightarrow}} u_0$, then 
$u_n(t) \underset{n \rightarrow \infty}{\overset{C^0}{\longrightarrow}} u(t)$.
Let $\eta>0$. \\
By $(\ref{maj_u_C1})$, there exists $C>0$
 such that for all $s\in (0,t]$, $\|u_x(s)\|_\infty \leq 
 \frac{C}{\sqrt{s}}$. \\
 So we can choose
  $\eta'>0$ such that $$\eta' e^{\int_0^t [\|u_x(s)\|_\infty^q+1] \;ds}\leq \eta.$$
 Let $n_0\geq 1$ such that for all $n\geq n_0$, $\|u_n-u_0\|_{\infty,[0,1]}\leq \eta' $.\\
Let $n\geq n_0$ and $s\in [0,t]$.
We denote $u_n(s)$ the solution at time $s$ of problem $(PDE_m)$ with initial condition $u_n$ and set :
$$z(s)=[u_n(s)-u(s)]e^{-\int_0^s [\|u_x(s')\|_\infty^q+1] \;ds'}.$$ 
 We see that $z$ satisfies
 \begin{equation}
 \label{equ_aux_z}
  z_s=x^{2-\frac{2}{N}}z_{xx}+b\, z_x+c\, z
 \end{equation}
where $$b(s,x)=u_n(s) \frac{(u_n)_x(s,x)^q-u_x(s,x)^q}{(u_n)_x(s,x)-u_x(s,x)} 
\mbox{ if }  (u_n)_x(s,x)\neq u_x(s,x) \mbox{ and $0$ else}$$ 
and  
 $$c=[{u_x}^q- \|u_x\|_\infty^q-1]<0.$$

  Since $z\in C([0,t]\times [0,1])$, $z$ reaches its maximum and its minimum. \\
  Assume that this maximum is greater than $\eta'$. Since $z=0$ for $x=0$ and $x=1$ 
and $z\leq \eta'$ for $s=0$, it can be reached only in 
  $(0,t]\times (0,1)$ but this is impossible because $c<0$ and (\ref{equ_aux_z}). 
  We make a similar reasoning for the minimum.
   Hence,
  $|z|\leq \eta'$ on  $[0,t]\times [0,1]$. \\Finally,
  $\|u_n-u\|_{\infty,[0,1]\times [0, t]} \leq \eta'
   e^{\int_0^t [\|(u_x(s))\|_\infty^q+1] \;ds}\leq \eta$ for all $n\geq n_0$. Whence the result.

\underline{Second step : } since $u_n \underset{n \rightarrow \infty}{\overset{C^1}{\longrightarrow}} u_0$, $\|u_n\|_{C^1}$ is bounded so there exists $K>0$ such that for all $n\geq 1$, $\mathcal{N}[u_n]\leq K$.
Then, from Lemma \ref{lem_compacite}, since $t>0$, $\{u_n(t),\;n\geq 1\}$ is relatively compact in $Y_m^1$ and has a single accumulation point $u(t)$ from first step. Whence the result.\hfill $\square$

\subsection{A strict Lyapunov functional for $(T(t))_{t\geq 0}$}

\subsubsection{Reminder on strict Lyapunov functionals}

We recall some definitions in the context of a continuous dynamical system $(T(t))_{t \geq 0}$ on $Y_m^1$, including strict Lyapunov functional and Lasalle's invariance principle.

\begin{Def}
 Let $u_0\in Y_m^1$.
\begin{itemize}
 \item  $\gamma_1(u_0)=\{T(t)u_0,\; t\geq 1 \}$ is the 
 trajectory of $u_0$ from $t=1$.
\item  $\omega(u_0)=\{v\in Y_m^1,\;
\exists t_n\rightarrow +\infty,\; t_n\geq 1,\; T(t_n)u_0 \underset{n \rightarrow +\infty}{\longrightarrow}v \mbox{ in }Y_m^1\}$\\ is the $\omega$-limit set of $u_0$.
\end{itemize}
\end{Def}

\begin{Def}\quad 
\begin{itemize}
\item [i)]$\mathcal{F}\in C(Y_m^1,\mathbb{R})$ is a Lyapunov functional if for all $u_0\in Y_m^1$,
$$t\mapsto \mathcal{F}[T(t)u_0] \mbox{ is nonincreasing on } [1,+\infty) .$$
\item [ii)]A Lyapunov functional $\mathcal{F}$ is a strict Lyapunov functional if 
$$ \mathcal{F}[T(t)u_0] =\mathcal{F}[u_0] \mbox{ for all } t\geq 0  \mbox{ implies that $u_0$ is an equilibrium point.}$$
\end{itemize}
\end{Def}

\begin{propo} Lasalle's invariance principle.\\
\label{equilibre}
 Let $u_0\in Y_m^1$.
Assume that the dynamical system $(T(t))_{t\geq 0}$ admits a strict Lyapunov functional
 and that $\gamma_1(u_0)$ is relatively compact in $Y_m^1$.\\
Then the $\omega$-limit set $\omega(u_0)$ is nonempty and consists of equilibria of the dynamical system.
\end{propo}

See \cite[p. 143]{Haraux} for a proof.

\subsubsection{Approximate Lyapunov functionals}
We recall that for all $\epsilon>0$ and $x>0$, 
$f_\epsilon(x)=(x+ \epsilon)^q-\epsilon^q$.
In order to introduce the approximate Lyapunov functional, we first need a double primitive $H_\epsilon$ of $\frac{1}{f_\epsilon}$ that converges uniformly to $ 
H(x)=\dfrac{x^{2-q}}{(2-q)(1-q)}$ on compacts of $\mathbb{R}^+$.  The next lemma provides it.

\begin{lem} 
Let $q\in(0,1)$ and $\epsilon \in(0,1]$. \\
For $x\geq 0$, we set 
$H_\epsilon(x)=\int_0^x\int_1^y \frac{dt}{f_\epsilon(t)}dy+\frac{x}{1-q}$ and $ 
H(x)=\dfrac{x^{2-q}}{(2-q)(1-q)}$.
\begin{itemize} 
 \item [i)]  $H_\epsilon$ is continuous on $[0;+\infty)$, twice differentiable on $(0;+\infty)$.\\
  $H_\epsilon''=\frac{1}{f_\epsilon}$ on $(0,+\infty)$.
 \item[ii)] $H_\epsilon$ converges uniformly to $ H$  on $[0;R]$ as $\epsilon$ tends to $0$, for
any $R>0$.
\end{itemize}
\end{lem}

\textit{Proof :} Let $R>0$ and $x\in [0,R]$. We begin with two remarks :\\
- Since $f_\epsilon\geq f_1$ and $f_1$ is concave, then, denoting $K=\frac{f_1(R)}{R}$, we have 
 $$f_\epsilon(t)\geq K t \mbox{\qquad for any } t>0.$$
- We also note that for any $t>0$, $$0\leq t^q-f_\epsilon(t)\leq \epsilon^q.$$
Indeed, setting $g(\epsilon)=(t+\epsilon)^q-\epsilon^q$, we have $$t^q-f_\epsilon(t)=g(0)-g(\epsilon)=-\epsilon\int_0^1 g'(\epsilon s) ds= \epsilon \int_0^1 \frac{q}{(s\epsilon)^{1-q}}-\frac{q}{(t+\epsilon s)^{1-q}} ds \leq \epsilon^q.$$

i) Let us set for $y>0 $, $\gamma_\epsilon(y)=\int_1^y 
\frac{dt}{f_\epsilon(t)}$. Since  for all $t>0$, 
$f_\epsilon(t)\geq K \, t$, then
$\gamma_\epsilon(y)=O(|\log(y)|)$. Hence $\gamma_\epsilon$ is integrable on $(0,R]$.
Then, $H_\epsilon$ is continuous on $[0,R]$. The other facts are obvious.\\
ii) We can write 
$$H_\epsilon(x)-H(x)=\int_0^x\int_1^y \frac{1}{f_\epsilon(t)}-\frac{1}{t^q} 
\,dt\,dy= \int_0^x\int_1^y \frac{t^q-f_\epsilon(t)}{t^q f_\epsilon(t)}dt dy.$$ Using the first two remarks, we get
$$|H_\epsilon(x)-H(x)|\leq \frac{\epsilon^q}{K}\int_0^R |\int_1^y \frac{dt}{t^{1+q}}|\leq \frac{\epsilon^q}{Kq}[\frac{R^{1-q}}{1-q}+R].$$
Whence the result.\hfill $\square$

\begin{Def}
 Let $\epsilon >0$. For $u\in Y_m^1$, we define 
$$ \mathcal{F}_\epsilon(u)=\int_0^1 H_\epsilon(u_x)-\dfrac{u^2}{2x^{2-q}} dx.$$
\end{Def}

We would like to remind to the reader that, if $u_0\in Y_m$ is given, $u^\epsilon$ denotes  the solution of problem
$(PDE_m^\epsilon)$ (see (\ref{equ_u_eps1})) with initial condition $u_0$.

\begin{lem} Let $u_0\in Y_m^1$. 
\label{lem_lyap1}
For all $0<t<s<T_{max}^\epsilon$, 
 $$ \mathcal{F}_\epsilon[u^\epsilon(s)]\leq \mathcal{F}_\epsilon[u^\epsilon(t)].$$
More precisely, for all $t>0$,
$$\frac{d}{dt} \mathcal{F}_\epsilon [u^\epsilon(t) ]=
-\int _0 ^1 \frac{(u^\epsilon)_t^2}{x^{2-q}f_\epsilon((u^\epsilon)_x) }dx. $$
\end{lem}

\textit{Proof :} Let $t>0$ and $\eta>0$ such that $I=[t-\eta,t+\eta]\subset (0,T_{max}^\epsilon)$.\\
By Theorem \ref{thm_existence_u_epsilon} iii), there exists $\mu >0$ such that $(u^\epsilon)_x\geq \mu $ on $I\times [0;1]$ so it will allow us to use that $H_\epsilon \in C^2([\mu , +\infty))$.\\
By Lemma \ref{lem_convergence_u_epsilon}, $(u^\epsilon)_{tx},\,(u^\epsilon)_{t},\,(u^\epsilon)_{x}$ are bounded on $I\times [0;1]$ and there exists $K>0$ such that  $|(u^\epsilon)_{xx}|\leq \frac{K}{x^{1-q}}$ on $I\times [0;1]$. 
Moreover,  $0\leq \frac{u^\epsilon}{x^{2-q}}\leq \frac{\|(u^\epsilon)_x\|_{C^1}}{x^{1-q}}\leq \frac{K'}{x^{1-q}}$  on $I\times [0;1]$ for some $K'>0$.
Note also that $u^\epsilon_{tx}=u^\epsilon_{xt}$ since $u^\epsilon\in C^2((0,T)\times (0,1])$.
All these facts allow us to differentiate $\mathcal{F}_\epsilon [u^\epsilon(t) ] $
and then to integrate by parts :
$$\frac{d}{dt} \mathcal{F}_\epsilon (u^\epsilon(t) )=
\int_0^1 (u^\epsilon)_{tx} H_\epsilon'((u^\epsilon)_x)-\frac{(u^\epsilon)_t u^\epsilon}{x^{2-q}} dx
$$ 
$$=\left[ (u^\epsilon)_t H_\epsilon'((u^\epsilon)_x)\right] _0^ 1-\int_0^1 (u^\epsilon)_t\left[ 
\frac{(u^\epsilon)_{xx}}{f_\epsilon((u^\epsilon)_x)}+\frac{ u^\epsilon}{x^{2-q}} \right]. $$
Hence, $\displaystyle \frac{d}{dt} \mathcal{F}_\epsilon (u^\epsilon(t) )=
-\int _0 ^1 \frac{(u^\epsilon)_t^2}{x^{2-q}f_\epsilon((u^\epsilon)_x) }$ since 
$(u^\epsilon)_t(t,0)=(u^\epsilon)_t(t,1)=0$. \hfill $\square$

\subsubsection{A strict Lyapunov functional}

\begin{Def}
  For $u\in Y_m^1$, we define 
$$ \mathcal{F}(u)=\int_0^1 H(u_x)-\dfrac{u^2}{2x^{2-q}} dx$$
where $H(v)=\dfrac{v^{2-q}}{(2-q)(1-q)}$ for all $v\in \mathbb{R}$.
\end{Def}

\begin{thm}\quad \\
\label{thm_lyap}
$\mathcal{F}$ is a strict Lyapunov functional for the dynamical system $T(t)_{t\geq 0}$ on $Y_m^1$.

\end{thm}

This fact cannot be obtained directly by the formal computation shown in the introduction since $u_x$ can vanish on a whole interval. Indeed, consider for instance an initial condition $u_0\in Y_M$ such that 
$ u_0\geq U_a$ where $a>A$. By comparison principle, for all $t\geq 0$, $M\geq u(t)\geq U_a$ so $u_x(t)=0$ at least on $[\frac{A}{a},1]$.  \\

\textit{Proof :} \underline{First,} by application of the dominated 
convergence theorem, since $H$ is continuous on 
$[0,\infty)$, it is clear that $\mathcal{F}$ is continuous 
on $Y_m^1$.\\
\underline{Secondly,} let $u_0\in Y_m^1$ and $0 < t \leq 
s $. Let us show that $\mathcal{F}(u(t))\geq 
\mathcal{F}(u(s))$.\\
Let $t'>0$. We first need to show that $\mathcal{F}_\epsilon(u^\epsilon(t')) \underset{\epsilon\rightarrow 0}
{\longrightarrow} \mathcal{F}(u(t'))$.\\
Lemma \ref{lem_convergence_u_epsilon} tells us that 
$u^\epsilon(t')\underset{\epsilon\rightarrow 0}{\longrightarrow} u(t')$ in $C^1([0,1])$.
In particular,
$(u^\epsilon(t'))_{\epsilon\in(0;1)}$ is bounded in $C^1([0;1])$ 
so we have a domination independent of $\epsilon$. Since $H_\epsilon$ converges 
uniformly
to $H$ on compact subsets of $[0;+\infty)$ and $H$ is continuous on $[0;+\infty)$, by the dominated 
convergence theorem,
we obtain easily that $\mathcal{F}_\epsilon(u^\epsilon(t')) \underset{\epsilon\rightarrow 0}
{\longrightarrow} \mathcal{F}(u_0(t'))$. \\
Then, as we know from Lemma \ref{lem_lyap1} that $ 
\mathcal{F}_\epsilon[u^\epsilon(t)]\geq 
\mathcal{F}_\epsilon[u^\epsilon(s)]$, the result follows by letting $\epsilon$ go to zero.\\
\underline{Thirdly,} denoting $R=\underset{t'\in[t;s]}{\sup} \|u_0(t') \|_{C^1}$, 
we want to show that 
\begin{equation}
\label{inequation_F}
 \mathcal{F}[u(t)]- \mathcal{F}[u_0(s)]\geq \frac{1}{(R+1)^q}\iint_{[t;s]\times 
[0;1]} (u_0)_t^2.
\end{equation}
By Lemma $\ref{lem_lyap1}$,
$$\mathcal{F}_\epsilon(u^\epsilon(t))- \mathcal{F}_\epsilon(u^\epsilon(s))=
\iint_{[t;s]\times 
[0;1]} \frac{(u^\epsilon)_t^2}{x^{2-q}f_\epsilon((u^\epsilon)_x)}.$$
By Lemma \ref{lem_convergence_u_epsilon}, $(u^\epsilon)_x $ tends to $u_x$ in $C([t,s]\times [0;1])$, so there exists $\epsilon_0>0$ such that 
for all $\epsilon\in(0;\epsilon_0)$, $\underset{t'\in[t;s]}{\sup} \|u^\epsilon(t') \|_{C^1}\leq R+1$. \\
Note that $\frac{1}{x^{2-q}}\geq 1$ for $0<x\leq 1$ and that 
$0<f_\epsilon(u^\epsilon_x)\leq (u^\epsilon_x)^q\leq 
(R+1)^q$ on $[t,s]\times [0,1]$. So, $$\mathcal{F}_\epsilon(u^\epsilon(t))- 
\mathcal{F}_\epsilon(u^\epsilon(s))\geq  
\frac{1}{(R+1)^q}\iint_{[t;s]\times 
[0;1]} (u^\epsilon_t)^2.$$
By Lemma \ref{lem_convergence_u_epsilon}, $(u^\epsilon)_t $ tends to 
$(u_0)_t$ in $C([t;s]\times [0,1])$, hence by taking the limit as $\epsilon$ 
goes 
to 0, we obtain the result.\\
\underline{Finally,} assume that $ \mathcal{F}[u(t)]= \mathcal{F}[u_0]$ for all $t\geq 0$.\\
Let $[t,T]\subset (0,\infty)$. Formula (\ref{inequation_F}) shows that 
$u_t=0$ on $[t,T]$. So $u_t=0$ on $(0,\infty)\times [0,1]$. Then, by continuity of $u$ on $[0,\infty)\times [0,1]$, we get 
that $T(t)u_0=u_0$ for all $t\geq 0$, i.e. $u_0$ is an equilibrium of the dynamical system $(T(t))_{t\geq 0}$.\hfill $\square$

\subsection{Convergence to a stationary state for $m\leq M$ : proof of Theorem \ref{thm_convergence}}

In the case of $m<M$, there is a unique stationary solution
for problem $(PDE_m)$ so the convergence  is not really surprising. \\
But for $m=M$ there is a continuum of stationary solutions (all $U_a|_{[0,1]}$ for $a\geq A$)
and the behaviour could be much more complicated. However, thanks to the good properties of the set of steady states (see Proposition \ref{prop_stat}) and since the problem is one-dimensional, convergence can  be shown by arguments in the spirit of \cite{Matano} or \cite{Zelenyak}.
\\

\textit{Proof of Theorem \ref{thm_convergence} :} Let $u_0\in Y_m$. Let us set $u_1=u_0(1)\in Y_m^1$.\\
To get the result, we just have to study $\underset{t\rightarrow +\infty}{\lim} T(t)u_1$. \\
Thanks to Lemma $\ref{lem_compacite}$, $\gamma_1(u_1)$ is relatively compact in $Y_m^1$ and 
since $\mathcal{F}$ is a strict Lyapunov functional for $(T(t))_{t\geq 0}$, we know by Lasalle's invariance principle (Proposition \ref{equilibre}) that 
the $\omega$-limit set $\omega(u_1) $ is non empty and contains only stationary 
solutions.\\
\underline{First case :} $m<M$. Then from Theorem $\ref{thm_stat}$, there exists a unique 
stationary solution $U_a$ with $a<M$. Hence, $\omega(u_1)=\{U_a\}$ so $T(t)u_1 \underset{t 
\rightarrow + \infty}{ \longrightarrow U_a}.$ \\
\underline{Second case :} $m=M$. Assume by contradiction that $\omega(u_1)$ contains two different stationary 
solutions $U_a$
and $U_b$ with $A\leq a<b$. Then we chose any $c\in(a,b)$. From Proposition  $\ref{prop_stat})$ ii), 
there exists $\epsilon>0$ such that $V_M(U_a,\epsilon)\prec U_c$ and $\{ u\in Y^1_M,\; u\leq U_c \}\cap 
V_M(U_b,\epsilon)=\emptyset$. \\
Since $U_a\in \omega(u_1)$, there exists $t_a$ such that $T(t_a)u_1 \in V_M(U_a,\epsilon)$. Hence, 
$T(t_a)u_1\leq U_c$ and then by comparison principle, for all 
$t\geq t_a$, $T(t)u_1\leq U_c$. But, since $U_b\in \omega(u_1)$, there exists $t_b\geq t_a$ such 
that $T(t_b)u_1\in V_M(U_b,\epsilon)$, and this is a contradiction because 
$\{ u\in Y^1_M,\; u\leq U_c \}\cap V_M(U_b,\epsilon)= \emptyset$. Hence, $\omega(u_1)$ is a singleton 
$\{U_a\}$
with $a\geq A$ so $T(t)u_1 \underset{t \rightarrow + \infty}{ \longrightarrow }U_a$. \hfill $\square$

\section{Finite time blow-up and self-similar solutions in supercritical case $m>M$}

In this section, we only consider the supercritical case, i.e. when $m>M$.\\
We shall prove that classical solutions of problem $(PDE_m)$ blow up in finite time. 
The idea of the proof is to exhibit a subsolution 
$\underline{u}(t,x)=V(a(t)x)$ which turns out to be a self-similar solution after some time. \\
This is why we are interested in the following ordinary differential equation.

\subsection{An auxiliary ordinary differential equation}

\begin{Def}
 Let $\epsilon>0$.\\
We define the problem $(Q_\epsilon)$ by :
\begin{eqnarray}
\label{Q_epsilon}
 x^{2-q} \ddot{V}+V \dot{V}^q&=&\epsilon x \dot{V} \mbox{ \qquad  } x>0  \\
\label{CI1}
 V(0)&=&0  \\
\label{CI2}
 \dot{V}(0)&=&1 
\end{eqnarray}

\end{Def}

\begin{Def} Let $\epsilon>0$ and $R>0$.\\
We say that $V$ is a solution of problem $(Q_{\epsilon})$ on $[0;R[$ 
if :
\begin{itemize}
\item $V\in C^1([0;R[)\cap C^2(]0;R[).$ 
\item $V$ is nondecreasing. 
\item $V$ satisfies \eqref{Q_epsilon} on $]0;R[$ and the conditions \eqref{CI1} and \eqref{CI2}.
\end{itemize}
This definition can obviously be adapted for the case of a closed interval $[0;R]$ or for 
$R=+\infty$.
\end{Def}

We summarize in the following theorem some very helpful results about solutions of problem $(Q_\epsilon)$.\\
Recall that $U_1$ is the solution of problem $(P_1)$ and that $A$ is the first point from which $U_1$ is constant.
 
 \begin{thm}
 \label{thm_V_epsilon}
 There exists a unique solution of problem $(Q_\epsilon)$ on $[0,\infty)$.\\
 If $\epsilon>0$ is small enough, $V_\epsilon$ is concave and there is a first point $A_\epsilon<\infty$ from which $V_\epsilon$ is 
 constant with value $M_\epsilon$ greater than $M$.\\
 Moreover : 
\begin{itemize}
\item [i)] $ \|V_\epsilon-U_1\|_{C^1([0,\infty))} \underset{\epsilon \rightarrow 0}
{\longrightarrow}0$.
\item [ii)]$A_\epsilon \underset{\epsilon \rightarrow 0}{\longrightarrow} A$.
\item [iii)]$V_\epsilon(A+1)\underset{\epsilon \rightarrow 0}{\longrightarrow}M$.\\
 That is to say that the constant reached by $V_\epsilon$ can be chosen as close to $M$ as we wish  provided that $\epsilon$ is small enough.
\end{itemize}
 \end{thm}

\textbf{Remark :} The fact that $M_\epsilon>M$ and that $V_\epsilon$ is concave for small $\epsilon>0$ follow from the proof of Theorem \ref{thm_self-similar}.\\

The proof of Theorem \ref{thm_V_epsilon} follows from the following successive lemmas.
We begin by giving an a priori property of the solutions of problem $(Q_\epsilon)$.

\begin{lem}
\label{lem_constant}
Let $\epsilon>0$ and $V$ a solution of problem $(Q_{\epsilon})$ on $[0;+\infty[$.\\
If for some $x_0\geq 0$, $\dot{V}(x_0)=0$, then for all $x\geq x_0$,  $V(x)=V(x_0)$. 
\end{lem}

\textit{Proof :} First, $x_0>0$ since $\dot{V}(0)=1$, then $V(x_0)>0$ since $V$ is nondecreasing.
Let $x_1=\sup \{ r\geq x_0, \; s.t.  \; V \mbox{ constant on } [x_0,r] \}$. \\
Assume by contradiction that $x_1<\infty$. 
Then, by continuity, $V(x_1)=V(x_0)>0$ and $\dot{V}(x_1)=0$. 
Writing equation (\ref{Q_epsilon}) as
$$ x^{2-q}\ddot{V}=\dot{V}^q(-V+\epsilon x \dot{V}^{1-q}),$$
we see that the first factor of the RHS is nonnegative and that the second one keeps negative for $x$ close enough to $x_1$ by 
continuity since $-V(x_1)+\epsilon x_1 \dot{V}(x_1)^{1-q}<0$. Hence 
$\ddot{V}$ is nonpositive for $x$ close enough to $x_1$. But $\dot{V}(x_1)=0$ and $\dot{V}\geq 0$, 
then $\dot{V}=0$ near of $x_1$, which contradicts the definition of $x_1$. So, $x_1=\infty$.\hfill $\square$\\

We now prove the local existence of a solution of problem $(Q_\epsilon)$.

\begin{lem}
\label{lem_existence_locale_V_eps}
 Let $\epsilon\in (0,1]$.\\
There exists   $\delta>0$ independent of $\epsilon$ such that the problem $(Q_\epsilon)$ admits a 
unique solution 
 on $[0,\delta]$.
\end{lem}

\textit{Proof :} The method used is a fixed point argument, as for the local existence of the 
solutions of problem $(P_a)$.\\
Let us define $$E=\{V\in C^1([0,\delta]),\; \dot{V}\geq 0,\; V(0)=0,\; \dot{V}(0)=1,\;\|\dot{V}-1 
\|_\infty\leq
 \dfrac{1}{2} \}.$$ 
$E$ equipped with the metric induced by the norm $\|V\|_{E}=\|\dot{V}\|_\infty$ is 
a complete metric space.
We define $F$ by :
 $$F:E\rightarrow C^1([0;\delta])$$  $$F(V)(x)=x-\int_0^x\int_0^y 
\dfrac{V(s)}{s}\dfrac{\dot{V}(s)^q}{s^{1-q}}ds\;
 dy+\epsilon \int_0^x\int_0^y \dfrac{\dot{V}}{s^{1-q}}ds\; dy.$$ 
Since for all $V\in E$, $\|\dot{V} \|_\infty \leq \frac{3}{2}$ and $0<\epsilon\leq 1$, we easily get 
that  $$\|F(V)'-1\|_\infty  \leq \left((\frac{3}{2})^{q+1} +\frac{3}{2}\right)\frac{\delta^q}
{q}\leq \frac{1}{2},$$ provided that $\delta$ is chosen small enough. 
Hence, $F$ sends $E $ into $E$. \\
We can apply the mean value theorem to function $z\mapsto z^q$, since for all $V\in E$, $\frac{1}
{2} \leq \dot{V}  \leq \frac{3}{2}$. Finally, we obtain for all $(V_1,V_2)\in E^2$ :
$$\| F(V_2)-F(V_1)\|_{E}
\leq \left(\dfrac{(\frac{3}{2}) ^q}{q}+\frac{\frac{3}{2}}{(\frac{1}{2})^{1-q}}+\frac{1}{q} 
\right)\delta^q \|V_2-V_1\| _{E} .
$$ 
Hence, if $\delta$ is small enough, $F$ is a contraction so there exists a fixed point $V$ of $F$. 
Since 
$F(V)\in C^2((0,\delta])$ when $V\in C^1([0;\tau])$ then $V$ is  a solution of 
$(Q_\epsilon)$.\\
Finally, it is easy to check that a solution of $(Q_\epsilon)$ is necessarily a fixed point of 
$F$, which proves the uniqueness. \hfill $\square$

\begin{lem} \label{lem_existence_presque_globale_V_epsilon}
 Let $\epsilon\in (0,1]$.\\
There exists a unique $A_\epsilon\geq \delta$ and  a unique maximal solution $V_\epsilon$ on $[0,A_\epsilon)$  of $(Q_\epsilon)$ such that 
$\dot{V_\epsilon}>0$ on $[0,A_\epsilon)$.
\end{lem}

\textit{Proof :} the proof follows from Lemma \ref{lem_existence_locale_V_eps}, similarly to that of Theorem \ref{thm_existence_globale_sol_stationnaires}.

\begin{lem}
 Let $\epsilon\in(0,1]$. 
\begin{itemize}
 \item [i)] We have the following formula :
for $x\in [0,A_\epsilon) $:\\
$$\dot{V_\epsilon}(x)=\exp(\frac{\epsilon(1-q) x^q}{q})
\left( 1-(1-q)\int_0^x \frac{V_\epsilon(s)}{s^{2-q}}\exp(-\frac{\epsilon(1-q)s^q}{q})ds  
\right)^{\frac{1}{1-q}} .$$
\item [ii)] For all $x\in [0,A_\epsilon) $, $$0<\dot{V_\epsilon}(x)\leq\exp(\frac{\epsilon(1-q) x^q}{q}),$$
$$0\leq V_\epsilon(x)\leq \int_0^x \exp(\frac{\epsilon(1-q) s^q}{q})ds.$$
\end{itemize}
\end{lem}

\textit{Proof :} 
$i)$ Let us set $w=\dot{V_\epsilon}$. On $(0,A_\epsilon)$, $w$ satisfies 
$$\dot{w}=-\frac{V_\epsilon(s)}{s^{2-q}}w^q+\frac{\epsilon}{x^{1-q}}w.$$
We recognize a Bernoulli type ordinary differential equation.
 Since $w>0$ on $[0,A_\epsilon)$, then we can divide by $w^q$.
Setting $z=w^{1-q}$, we obtain $$\dot{z}=\frac{\epsilon(1-q)}{x^{1-q}}z-\frac{V_\epsilon(s)(1-q)}{s^{2-q}}$$
which can be easily integrated. Whence the formula.\\
$ii)$ $V_\epsilon$ is increasing and $V_\epsilon(0)=0$ so $V_\epsilon\geq 0$ on $[0,A_\epsilon)$ 
whence the results. \hfill $\square$

\begin{rmq}
\label{rmq_negatif}
A consequence of last lemma is that \\if $1-(1-q)\int_0^x \frac{V_\epsilon(s)}{s^{2-
q}}\exp(-\frac{\epsilon(1-q)s^q}{q})ds \leq 0$, then $A_\epsilon \leq x$.
\end{rmq}

\begin{lem}
Let $\epsilon\in(0,1]$. \\
There exists a unique solution $V_\epsilon$ of problem $(Q_\epsilon)$ defined on $[0,+\infty)$.
\end{lem}

\textit{Proof :} If $A_\epsilon=\infty$ then all is already done. \\
Else, if $A_\epsilon<\infty$ then  
$\dot{V_\epsilon}(x)$ has to go to zero as $x$ goes to $A_\epsilon$. Indeed, $(V_\epsilon(x),
\dot{V_\epsilon}(x))$ must go out of any compact of $\mathbb{R}\times (0,+\infty)$ as $x$ goes to 
$A_\epsilon$ but by $ii)$ of the above lemma $V_\epsilon(x)$ and $\dot{V_\epsilon}(x)$ keep bounded 
for $x$ bounded. So $\dot{V_\epsilon}(x)\underset{x \rightarrow A_\epsilon}{\longrightarrow}0$.\\
Hence, by Cauchy criterion, $V_\epsilon(x)$ has a limit $L_\epsilon$ as $x$ goes to $A_\epsilon$. 
Moreover, by the equation $(\ref{Q_epsilon})$, $\ddot{V_\epsilon}(x)$ goes to zero as $x$ goes to 
$A_\epsilon$. Since the constants are solutions of $(\ref{Q_epsilon})$, then $V_\epsilon$ can be 
extended by the constant $L_\epsilon$ on $[A_\epsilon,+\infty)$ to a $C^2$ function on $(0,
+\infty)$ which is a solution of problem $(\ref{Q_epsilon})$. \\
For proving the uniqueness, let $V$ a solution on $[0,+\infty)$. By uniqueness of the solution on 
$[0,A_\epsilon)$ then $V=V_\epsilon$ on $[0,A_\epsilon)$. Then, by continuity, 
$V(A_\epsilon)=V_\epsilon(A_\epsilon)$ and 
$\dot{V(}A_\epsilon)=0$. But now by Lemma \ref{lem_constant}, $V$ is constant on $[A_\epsilon,
\infty)$ so $V=V_\epsilon$ on $[0,\infty)$. \hfill $\square$

\begin{lem}\quad
\label{lem_convergence_V_eps}
\begin{itemize}
\item [i)] $A_\epsilon<\infty$ for $\epsilon$ small enough. \\
Moreover $A_\epsilon \underset{\epsilon \rightarrow 0}{\longrightarrow}A$. 
\item [ii)]$ \|V_\epsilon-U_1\|_{C^1([0,\infty))} \underset{\epsilon \rightarrow 0}
{\longrightarrow}0$.\\
As a consequence, $V_\epsilon(A+1)\underset{\epsilon \rightarrow 0}{\longrightarrow}M$.
\end{itemize}
\end{lem}

\textit{Proof :}\\
 \underline{First step :} We show that there exists $\gamma\in(0,\delta]$ independent of $\epsilon$ 
 such that 
$$ \| V_\epsilon-U_1\|_{C^1([0,\gamma])}\underset{\epsilon \rightarrow 0}{\longrightarrow}0.$$ 
 Let us set $\gamma=\min(\delta,\delta',A )$ where $\delta$ is a short existence time for all 
 $V_\epsilon$ and $\delta'=(\frac{q}{2})^{\frac{1}{q}}$.
 $\dot{U_1}>0$ on $[0,A)$ so equation $(\ref{ode})$ can be written on $[0,A)$ as 
 $$ \frac{d}{dx}\frac{\dot{U_1}^{1-q}}{1-q}=-\frac{U_1}{x^{2-q}}.$$
Hence, we have the following formula 
\begin{equation}
\label{formule_U_1}
\dot{U_1}(x)= \left(1-(1-q)\int_0^x \frac{U_1(s)}{s^{2-q}}ds\right)^{\frac{1}{1-q}},
\end{equation}
which is valid for all $x\in[0,A]$, by continuity.\\
Let $x\in [0,\gamma]$. 
$$\dot{V_\epsilon}(x)-\dot{U_1}(x)=$$ $$\left(e^{\epsilon\frac{(1-q)x^q}{q}}-1\right)
\left( 1-(1-q)\int_0^x \frac{V_\epsilon(s)}{s^{2-q}}\exp(-\frac{\epsilon(1-q)s^q}{q})ds  
\right)^{\frac{1}{1-q}}$$
$$+\left( 1-(1-q)\int_0^x \frac{V_\epsilon(s)}{s^{2-q}}\exp(-\frac{\epsilon(1-q)s^q}{q})ds  
\right)^{\frac{1}{1-q}}$$ 
$$- \left(1-(1-q)\int_0^x \frac{U_1(s)}{s^{2-q}}ds\right)^{\frac{1}{1-q}}.$$ 
Then,
$$|\dot{V_\epsilon}(x)-\dot{U_1}(x)|\leq (e^{\epsilon\frac{(1-q)x^q}{q}}-1)+\int_0^x \Big| 
\frac{V_\epsilon(s)}{s^{2-q}} e^{-\frac{\epsilon(1-q)s^q}{q}} 
- \frac{U_1(s)}{s^{2-q}}\Big|ds$$
 because $|a^{\frac{1}{1-q}}-b^{\frac{1}{1-q}}|\leq \frac{1}{1-q}| a-b|$ for all $(a,b)\in[0,1]^2$  
 since $\frac{1}{1-q}>1$.\\
 $$|\dot{V_\epsilon}(x)-\dot{U_1}(x)|\leq (e^{\epsilon\frac{(1-q)\gamma^q}{q}}-1)+\int_0^x  \frac{|
 V_\epsilon(s)-U_1(s)|}{s^{2-q}} e^{-\frac{\epsilon(1-q)s^q}{q}}ds$$ 
$$+ \int_0^x\frac{U_1(s)}{s^{2-q}}(1-e^{-\frac{\epsilon(1-q)s^q}{q}})ds.$$
Then, we get
$$\|\dot{V_\epsilon}-\dot{U_1}\|_{\infty,[0,\gamma]} \leq (e^{\epsilon\frac{(1-q)\gamma ^q}{q}}-1)$$ 
$$+\|
\dot{V_\epsilon}-\dot{U_1}\|_{\infty,[0,\gamma]} \frac{\gamma ^q}{q} 
 +(1-e^{-\frac{\epsilon(1-q)\gamma^q}{q}}) \int_0^\gamma \frac{U_1(s)}{s^{2-q}}ds.$$
Since $\frac{\gamma^q}{q}\leq \frac{1}{2}$ and $\int_0^\gamma \frac{U_1(s)}{s^{2-q}}ds\leq \frac{1}
{1-q}$, then
$$\|\dot{V_\epsilon}-\dot{U_1}\|_{\infty,[0,\gamma]} \leq 2\left[e^{\epsilon\frac{(1-q)\gamma^q}{q}}-1
+\frac{1-e^{-\frac{\epsilon(1-q)\gamma^q}{q}}}{1-q}\right]  . $$
Hence, $ \| V_\epsilon-U_1\|_{C^1([0,\gamma])}\underset{\epsilon \rightarrow 0}{\longrightarrow}0$.
 
 \underline{Second step :}
 Let $A'<A$. \\
 Let us show that for $\epsilon$ small enough $A_\epsilon\geq A'$ and that 
$$\|V_\epsilon-U_1\|_{C^1([0,A'])} \underset{\epsilon \rightarrow 0}{\longrightarrow}0.$$
Let us denote $V_{a,b}$ the solution of $(\ref{Q_epsilon})$ such that 
$V_{a,b}(\gamma)=a$ and $\dot{V}_{a,b}(\gamma)=b$.
 $U_1$ is the solution of equation $(\ref{Q_epsilon})$ for $\epsilon=0$ and initial condition 
 $(a,b)=(U_1(\gamma),\dot{U_1}(\gamma))$. 
Since  $\dot{U_1}>0$ on $[\gamma,A']$, the classical Cauchy-Lipschitz theory is here available, and 
by continuity of the solutions on $[\gamma,A']$ with respect to the parameter $\epsilon$ and the 
initial condition $(a,b)$, we know that if $\epsilon$ is small enough and if the initial condition 
$(a,b)$ is close enough to $(U_1(\gamma),\dot{U_1}(\gamma))$ then $\|V_{a,b}-U_1\|
_{C^1([\gamma,A'])}$ is as small as  we wish. But, thanks to the first step, taking $\epsilon$ even 
smaller, $(V_\epsilon(\gamma),\dot{V_\epsilon}(\gamma))$ can be made as close as necessary of 
$(U_1(\gamma),\dot{U_1}(\gamma))$. Finally, 
$\|V_\epsilon-U_1\|_{C^1([0,A'])}$ is as small as necessary when $\epsilon$ is small enough. This 
implies of course that $A_\epsilon\geq A'$ for $\epsilon$ small enough. Whence the results.

\underline{Third step :} Let us show that $ \|V_\epsilon-U_1\|_{C^1([0,A+1])} \underset{\epsilon 
\rightarrow 0}{\longrightarrow}0$.\\
Let $\alpha>0$ and $\eta=(\frac{\alpha}{4})^{1-q}$. \\
By formula $(\ref{formule_U_1})$, we have $(1-q)\int_0^A \frac{U_1(s)}{s^{2-q}}ds=1$ and since $\dot{U_1}=0$ on $[A,\infty)$, there 
exists $A'<A$ such that $1-(1-q)\int_0^{A'} \frac{U_1(s)}{s^{2-q}}ds \leq \eta$ and $ \| 
\dot{U_1}\|_{\infty,[A',\infty]}\leq \frac{\alpha}{2}$. \\
Let $x\in [A',\min(A_\epsilon,A+1)]$. We know that \\$\dot{V_\epsilon}(x)=\exp(\frac{\epsilon(1-q) x^q}
{q})
\left( 1-(1-q)\int_0^x \frac{V_\epsilon(s)}{s^{2-q}}\exp(-\frac{\epsilon(1-q)s^q}{q})ds  
\right)^{\frac{1}{1-q}}$\\ then $\dot{V_\epsilon}(x)\leq \exp(\frac{\epsilon (1-q)(A+1)^q}{q})
\left( 1-(1-q)\int_0^{A'} \frac{V_\epsilon(s)}{s^{2-q}}ds  \right)^{\frac{1}{1-q}} $ since 
$V_\epsilon\geq 0$. \\
There exists $\epsilon_1>0$ such that if $0<\epsilon\leq \epsilon_1$, $\exp(\frac{\epsilon (1-q)(A+1)^q}
{q})\leq 2$, so
$\dot{V_\epsilon}(x)\leq 2 \eta^\frac{1}{1-q}=\frac{\alpha}{2}$. Then $\|\dot{V_\epsilon}\|
_{\infty,[A',A+1]}\leq \frac{\alpha}{2}$ since $\dot{V_\epsilon}=0$ in $[A_\epsilon,\infty]$. 
Finally, if 
$0<\epsilon\leq \epsilon_1$, then $\|\dot{V_\epsilon}-\dot{U_1}\|_{\infty,[A',A+1]}\leq 
\frac{\alpha}{2}+\frac{\alpha}{2}=\alpha$. From the second step, let $\epsilon_2>0$ such that $ \|
\dot{V_\epsilon}-\dot{U_1}\|_{\infty,[0,A']}\leq \alpha$ for $0<\epsilon\leq \epsilon_2$. Let us 
set $\epsilon_0=\min(\epsilon_1,\epsilon_2)$. Then, if  if 
$0<\epsilon\leq \epsilon_0$, then $\|\dot{V_\epsilon}-\dot{U_1}\|_{\infty,[0,A+1]}\leq \alpha$.
\\
\underline{Fourth step :} Let us show that for $\epsilon$ small enough, $A_\epsilon<\infty$ and even  that $A_\epsilon 
\underset{\epsilon \rightarrow 0}{\longrightarrow}A$.\\
Let $A'\in(A,A+1]$. 
Since, by formula $(\ref{formule_U_1})$, we have $$1-(1-q)\int_0^{A'} \frac{U_1(s)}{s^{2-q}}ds=-(1-q)\int_A^{A'}\frac{U_1(s)}{s^{2-q}}<0  $$ 
and since
$$ \|V_\epsilon-U_1\|_{C^1([0,A'])} \underset{\epsilon \rightarrow 0}{\longrightarrow}0,$$
 then there exists $\epsilon_0>0$ such that for $0<\epsilon<\epsilon_0$, 
$$1-(1-q)\int_0^{A'} \frac{V_\epsilon(s)}{s^{2-q}}\exp(-\frac{\epsilon(1-q)s^q}{q})ds<0.$$
Let $0<\epsilon<\epsilon_0$. 
From remark $\ref{rmq_negatif}$, then $A_\epsilon\leq A'<\infty$.
Combined with step 2, this proves that $A_\epsilon \underset{\epsilon \rightarrow 0}
{\longrightarrow}A$\\ 
\underline{Last step :} Let $\epsilon>0$ small enough so that $A_\epsilon\leq A+1$. Then, 
 $\dot{V_\epsilon}=\dot{U_1}=0$ on $[A+1,\infty)$ thanks to Lemma $\ref{lem_constant}$. Moreover,
  thanks to step 3, $V_\epsilon(A+1)\underset{\epsilon \rightarrow 0}{\longrightarrow}U_1(A+1)=M$.
  Whence $ \|V_\epsilon-U_1\|_{C^1([0,\infty))} \underset{\epsilon 
  \rightarrow 0}  {\longrightarrow}0$. \hfill $\square$

\subsection{Blow-up and existence of self-similar solutions : proof of Theorems \ref{blow-up} and \ref{thm_self-similar}}

\underline{\textit{Proof of Theorem \ref{blow-up} :}} Let $m>M$ and $u_0\in Y_m$. \\
From Theorem \ref{thm_V_epsilon}, we can fix $\epsilon>0$ such that $V_\epsilon$ is constant equal to $L$ 
on $[A_\epsilon,\infty)$ with $L<m$.\\
Pick $t_0\in(0,T_{max}(u_0))$. From Theorem $\ref{thm_existence_u}$ iii), we know that $u_x(t_0,0)>0$, then $\alpha=\underset{x\in(0,1]}{\min} \dfrac{u(t_0,x)}{x}>0$. We now set $$a_0=\dfrac{\alpha}{\|
\dot{V_\epsilon}\|_{\infty,[0,+\infty)}}>0.$$
Then for all $x\in [0,1]$, 
$$V_\epsilon(a_0\,x )\leq 
\|\dot{V_\epsilon}\|_{\infty,[0,+\infty)}\,a_0\, x\leq 
\frac{a_0\,\|\dot{V_\epsilon}\|_{\infty,[0,+\infty)}}{\alpha}u(t_0,x).$$

Hence, $a_0$ is small enough so that
$$ V_\epsilon(a_0 x)\leq u(t_0,x)\; \mbox{ for all } x\in[0,1].$$
Let us set $$a(t)=\dfrac{a_0}{(1-\epsilon\, a_0^q\, q t)^{\frac{1}{q}}}$$ and $$T^*= \dfrac{1}
{\epsilon  \;a_0^q\; q}.$$
 Remark that $a(0)=a_0$ and
 $$\dot{a}(t)=\epsilon \; a(t)^{1+q} \mbox{ for } t\in [0,T^*).$$  
 We shall show that $$\underline{u}(t,x)=V_\epsilon(a(t)x)$$
  is a subsolution for $(PDE)_m$ with initial condition $u_0(t_0)$.\\ Let us set $T=\min(T_{max}
  (u_0)-t_0,T^*)$.\\
 Indeed, $\underline{u}(0,x)=V_\epsilon(a_0 x)\leq u_0(t_0,x)$ for all $x\in[0,1]$. 
 For all $t\in [0,T)$, we see that 
 $\underline{u}(t,0)=0= u_0(t_0+t,0) $ and
$\underline{u}(t,1)\leq L\leq m= u_0(t_0+t,1) $.\\
Moreover,  a straightforward calculation shows that 
$$\underline{u}_t-x^{2-q}\underline{u}_{xx}-\underline{u}\;\underline{u}_x\,^q=a^q\left[\frac{\dot{a}}{a^{1+q}} y\dot{V_\epsilon}(y)-
y^{2-q}\ddot{V_\epsilon}(y)-V_\epsilon(y) \dot{V_\epsilon}(y)^q \right]=0$$
where $y=a(t)x$.\\

From the comparison principle (cf Lemma $\ref{PC_PDE}$), $\underline{u}(t)\leq u_0(t_0+t)$ for all $t\in [0,T)$.\\
Now, if we assume that $T_{max}(u_0)>t_0+T^*$, then $T=T^*$ and by letting $t$ go to $T^*$, since 
$a(t)\underset{t \rightarrow T^*}{\longrightarrow} +\infty$, we obtain
$$L \leq u_0(t_0+T^*,x) \mbox{ for all } x\in(0,1].$$
Since $u_0(t_0+T^*,0)=0$, this contradicts the continuity of $u_0(t_0+T^*)$ at $x=0$. Hence, $T_{max}(u_0)\leq t_0+T^*<\infty$. \hfill $\square$ \\

\underline{\textit{Proof of Theorem \ref{thm_self-similar} :}} 
 From Theorem \ref{thm_V_epsilon}, we know the existence of 
$\epsilon_1>0$ such that for $\epsilon\in (0,\epsilon_1]$, $A_\epsilon\leq A+1$ and $V_\epsilon$ is flat from $x=A_\epsilon$.
We set $u_{\epsilon,a_0}(t,x)=V_\epsilon(a(t)x)$   where $ a(t)=\dfrac{a_0}{(1-\epsilon\, a_0^q\, q t)^{\frac{1}{q}}} $ and $a_0\geq A_\epsilon$.
The calculation in the proof of the previous theorem and the fact that 
$V_\epsilon$ is constant from the point $x=A_\epsilon$ prove that $\left(u_{\epsilon,a_0}\right)_{a_0\geq A_\epsilon}$ is a family of  solutions of problem $(PDE_{M_\epsilon})$ where $M_\epsilon=V_\epsilon(A+1)$.\\
 Now, by using the same methods as in the proof of Lemma 
\ref{lem_convergence_V_eps}, we can prove that 
\begin{equation}
\label{convergence_C1_V_eps}
\epsilon \mapsto V_\epsilon \mbox{ is continuous from } [0,\epsilon_1] \mbox{ to }  C^1([0,+\infty)).
\end{equation}
Then, $\epsilon\mapsto M_\epsilon=V_\epsilon(A+1)$ is continuous in $[0,\epsilon_1]$ so its image is a compact interval $I$. Since $u_{\epsilon,a_0}$  blows up, we  necessarily
have  $M_\epsilon>M$ if $\epsilon\in (0,\epsilon_1]$
so $I=[M,M^+]$ with $M^+>0$.\\
Finally, denoting $K_\epsilon=\underset {x\in(0,\infty)}{\sup} \frac{V_\epsilon(x)}{x}$, it is clear that $K_\epsilon=\underset {x\in(0,A_\epsilon]}{\sup} \frac{V_\epsilon(x)}{x}$ since $V_\epsilon$ is flat from $x=A_\epsilon$. Then, we have $$\mathcal{N}[u_{\epsilon,a_0}(t)]=K_\epsilon \,a(t)=\frac{K_\epsilon}{q^{\frac{1}{q}}\epsilon^{\frac{1}{q}}
(T_{max}- t)^{\frac{1}{q}}}, $$ where $$T_{max}=\frac{1}{q\epsilon a_0^q}.$$
We now prove that $V_\epsilon$ is concave for small $\epsilon$, which implies that $K_\epsilon=1$ so that the blow-up speed is known explicitly.
Let $\epsilon\in [0,\epsilon_1]$.\\ By Lemma \ref{lem_existence_locale_V_eps}, 
there exists $\delta>0$ independent of $\epsilon$ such that 
$\dot{V_\epsilon}\geq \frac{1}{2}$ on $[0,\delta]$. Since  $V_\epsilon$ is 
moreover nondecreasing on $[0,A+1]$, 
then $\underset{x\in(0,A+1]}{\min} \frac{V_\epsilon}{x}\geq \frac{\delta}{2}$.
By (\ref{convergence_C1_V_eps}), there exists $C>0$ such that $\|
\dot{V_\epsilon}\|_{\infty}\leq C$ for all $\epsilon \in [0,\epsilon_1]$.\\
Let $\epsilon_2=\min(\epsilon_1, \frac{\delta}{2\,C^{1-q}})$, $\epsilon \in(0,
\epsilon_2]$ and $x\in[0,A+1]$. \\We have $\epsilon\, x\,\dot{V_\epsilon}(x)^{1-q}\leq 
\epsilon_2\,x\, C^{1-q}\leq \frac{\delta}{2}x\leq V_\epsilon(x)$ and since 
$$x^{2-q}\ddot{V}_\epsilon(x)=-\dot{V_\epsilon}
(x)^q\left[V_\epsilon(x)-\epsilon x 
\dot{V_\epsilon}(x)^{1-q}\right],$$
we deduce that $V_\epsilon$ is concave and $K_\epsilon=\dot{V_\epsilon}(0)=1$ for $\epsilon\in [0,\epsilon_2]$. \hfill $\square$\\

{\bf Acknowledgements :} the author would like to thank Philippe Souplet 
for all discussions, encouragements and comments about this paper and also the referees whose remarks and suggestions helped a lot to improve this paper.


\begin{thebibliography}{00}
\bibitem{BCD}P. Biler, L. Corrias, J. Dolbeault, Large mass self-similar 
solutions of the parabolic-parabolic Keller-Segel model of chemotaxis. J. Math. 
Biol. 63 (2011), no. 1, pp. 1-32
\bibitem{BHN} P. Biler, D. Hilhorst and T. Nadzieja, Existence and nonexistence of solutions for a model of gravitational interaction of particles. II. Colloq. Math. 67 (1994), no. 2, 297-308
\bibitem{BKLN} P. Biler, G. Karch, P. Lauren\c cot and T. Nadzieja, The $8\pi$-problem for radially symmetric solutions of a chemotaxis model in a disk, Topol. Methods Nonlinear Anal. 27 (2006), no. 1, pp. 133-147
\bibitem{BKLN2} P. Biler, G. Karch, P. Lauren\c cot, and T. Nadzieja, The 8$\pi$-problem for radially symmetric                  
 solutions of a chemotaxis model in the plane, Math. Methods Appl. Sci., 29 (2006), pp. 1563-1583.
\bibitem{BCL} A. Blanchet, J.A. Carrillo and P. Lauren\c cot, Critical mass for a Patlak-Keller-Segel model with degenerate diffusion in higher dimensions, Calc. Var. Partial differential equations(2009) 35, no. 2, pp. 133-168
\bibitem{Blanchet2} A. Blanchet, J.A. Carrillo, and N. Masmoudi, Infinite time aggregation for the critical
    two-dimensional Patlak-Keller-Segel model, Comm. Pure Appl. Math., 61 (2008), pp. 1449-1481.
\bibitem{Blanchet3}A. Blanchet, J. Dolbeault, and B. Perthame, Two dimensional Keller-Segel model: Opti-
    mal critical mass and qualitative properties of solutions, Electron. J. Differential Equations,
    44 (2006), pp. 1-32.
    \bibitem{BL} A. Blanchet, P. Lauren\c cot, Finite mass self-similar blowing-up solutions of a chemotaxis system with non-linear diffusion, (2012) Communications in Pure and Applied Analysis, Vol. 11, Number 1
\bibitem{BRB} J. Bedrossian, N. Rodriguez and A.L Bertozzi,
Local and global well-posedness for aggregation equations and Patlak-Keller-Segel models with degenerate diffusion, Nonlinearity 24 (2011) 1683-1714 
\bibitem{Chavanis} P.H. Chavanis, Nonlinear mean field Fokker-Planck equations. Application  to the chemotaxis of biological populations, Eur. Phys. J. B 62, 2008,  pp. 179-208

\bibitem{CALCOR} V. Calvez, L. Corrias, The parabolic-parabolic Keller-Segel model in $\mathbb{R}^2$. Commun. MAth. Sci. 6(2008), no.2, 417-447

\bibitem{Cieslak_Winkler}T. Cieslak and M. Winkler,  Finite-time blow-up in a quasilinear system of chemotaxis. Nonlinearity 21 (2008), no. 5, 1057--1076.

\bibitem{DW}
K. Djie, M. Winkler, 
Boundedness and finite-time collapse in a chemotaxis system with volume-filling effect, Nonlinear Anal. TMA 72 (2010), pp. 1044-1064.
\bibitem{Dolbeault} J. Dolbeault and B. Perthame, Optimal critical mass in the two-dimensional Keller-Segel
    model in R2 , C. R. Math. Acad. Sci. Paris, 339 (2004), pp. 611-616.

\bibitem{Friedman} A. Friedman, Partial differential equations of parabolic type. Prentice-Hall, Inc., Englewood Cliffs, N.J. 1964 xiv+347 pp.
\bibitem{Haraux} T. Cazenave and A. Haraux, An introduction to semilinear evolution equations. Oxford Lecture Series in Mathematics and its Applications, 13. The Clarendon Press, Oxford University Press, New York, 1998. xiv+186 pp.
\bibitem{Herrero} M.A. Herrero, The mathematics of chemotaxis, Handbook of differential equations: evolutionary equations. Vol. III, 137-193, 
Handb. Differ. Equ., Elsevier/North-Holland, Amsterdam, 2007. 
\bibitem{Herrero-Sastre} M.A. Herrero and L. Sastre, Models of aggregation in dictyostelium discoideum : on the track of spiral waves, Networks and heterogeneous media, Volume 1, Number 2, June 2006, pp. 241-258
\bibitem{Herrero-Velazquez} M.A. Herrero and J.L. Velazquez, Singularity patterns in a chemotaxis model, Math. Ann. 306, pp. 583-623 (1996)
\bibitem{Horstmann1} D. Horstmann, From 1970 until present : the Keller-Segel model in chemotaxis and its consequences I, Jahresber. Deutsch. Math.-Verein., 105 (2003), pp. 103-165
\bibitem{Horstmann2} D. Horstmann, From 1970 until present : the Keller-Segel model in chemotaxis and its consequences II, Jahresber. Deutsch. Math.-Verein., 106 (2004), pp. 51-69
\bibitem{Horstmann-Winkler} D. Horstmann and M. Winkler, Boundedness vs. blow-up in a chemotaxis system, Journal of Differential Equations 215 (2005) pp. 52-107
\bibitem{HP} T. Hillen, K. J. Painter A user's guide to PDE models for chemotaxis. J. Math. Biol. 58 (2009), no. 1-2, 183--217. 
\bibitem{Hu} B. Hu, Blow-up theories for semilinear parabolic equations. 
Lecture Notes in Mathematics, 2018. Springer, Heidelberg, 2011. x+125 pp.
\bibitem{Kavallaris-Souplet} N.I. Kavallaris and P. Souplet, Grow-up rate and refined asymptotics for a two-dimensional Patlak-Keller-Segel model in a disk. 
SIAM J. Math. Anal. 40 (2008/09), no. 5, 1852-1881. 

\bibitem{Kim_Yao} I. Kim and Y. Yao, The Patlak-Keller-Segel model and its variations: properties of solutions via maximum principle. SIAM J. Math. Anal. 44 (2012), no. 2, 568--602.

\bibitem{KS} E.F. Keller and L.A. Segel, Initiation of slime mold aggregation viewed as an instability, J. Theor. Biol., 26 (1970), pp. 399-415
\bibitem{LS}
S. Luckhaus, Y. Sugiyama,
Asymptotic profile with the optimal convergence rate for a parabolic equation of chemotaxis in super-critical cases,
Indiana Univ. Math. J. 56 (2007) 1279-1297.
\bibitem{Lu} A. Lunardi, Analytic semigroups and optimal regularity in parabolic problems. Progress in Nonlinear Differential Equations and their Applications, 16. Birkh\"{a}user Verlag, Basel, 1995. xviii+424 pp.
\bibitem{Matano} H. Matano, Convergence of solutions of one-dimensional 
semilinear parabolic equations, J. Math. Kyoto Univ. 18 (1978), no. 2, pp. 221-227
\bibitem{MS}
N. Mizoguchi, Ph. Souplet,
Nondegeneracy of blow-up points for the parabolic Keller-Segel System,
preprint (2012).
\bibitem{Montaru} A. Montaru, Radially symmetric solutions of a semilinear parabolic-elliptic Patlak-Keller-Segel system with nonlinear sensitivity, \textit{to appear.}
\bibitem{Nagai}
T. Nagai, 
Blow-up of radially symmetric solutions to a chemotaxis system,
Adv. Math. Sci. Appl. 5 (1995), 581-601.
\bibitem{Pa} A. Pazy, Semigroups of linear operators and applications to partial differential equations. 
Applied Mathematical Sciences, 44. Springer-Verlag, New York, 1983. viii+279 pp.
\bibitem{Perthame} B. Perthame, PDE models for chemotactic movements: Parabolic, hyperbolic and kinetic,
Appl. Math., 49 (2004), pp. 539-564.
\bibitem{PW} M.H. Protter and H.F. Weinberger, Maximum principles in differential equations.  Springer-Verlag, New York, 1984. x+261 pp.
\bibitem{QS} P. Quittner and P. Souplet, Superlinear parabolic problems. Blow-up, global existence and steady states. Birkh\"{a}user Advanced Texts.  
 Birkh\"{a}user Verlag, Basel, 2007.
\bibitem{Patlak} C.S. Patlak, Random walk with persistence and external bias, Bull. Math. Biol. Biophys.,15 (1953), pp. 311-338  
\bibitem{SK}
Y. Sugiyama, H. Kunii, 
Global existence and decay properties for a degenerate Keller-Segel model with a power factor in drift term, 
J. Differential Equations 227 (2006), 333-364.
\bibitem{Suzuki}
T. Suzuki, Free energy and self-interacting particles. Progress in Nonlinear Differential Equations and their Applications, 62. Birkh\"{a}user Boston, Inc., Boston, MA, 2005.
\bibitem{Yao_Bertozzi} Y. Yao and  A.L. Bertozzi, Blow-up dynamics for the aggregation equation with degenerate diffusion, preprint (2012).
\bibitem{Zelenyak} T.I. Zelenyak, Stabilization of solutions of boundary value 
problems for a second order parabolic equation with one space variable (1968), 
Differential Equations, 4 , pp. 17-22
\end{thebibliography}
 \end{document}

 \section{Appendix}
 In order to lighten the reading, we chose to leave aside the (standard) proofs of two results. For the sake of completeness, we present them now.

\subsection{Proof of Lemma \ref{lem_existence_presque_globale_V_epsilon}}

Lemma \ref{lem_existence_locale_V_eps} gives a unique classical 
solution $V_0$ of $(Q_\epsilon)$ on $[0;\delta]$. We can now consider the following ordinary 
differential equation on the interval $[\dfrac{\delta}{2},+\infty[$ :
\begin{eqnarray} 
\ddot{V}&=&-\dfrac{V\dot{V}^q}{x^{2-q}}+\frac{\epsilon}{x^{1-q}}\dot{V} \nonumber \\ 
 V(\dfrac{\delta}{2})&=&V_0(\dfrac{\delta}{2}) \nonumber \\
 \dot{V}(\dfrac{\delta}{2})&=&\dot{V_0}(\dfrac{\delta}{2})>0 \nonumber
\end{eqnarray}
Setting $Z=(V,\dot{V})$, it becomes a first order nonautonomous ordinary differential equation 
$$\dot{Z}=F(x,Z,\epsilon)$$  $$\mbox{where }F(x,V,W,\epsilon)=\left(W,
\dfrac{-VW^q}{x^{2-q}}+\frac{\epsilon}{x^{1-q}}W\right).$$\\
Let us denote $\Omega=\mathbb{R} \times (0;+\infty)$. Since $F$ is locally 
Lipschitz with respect to $Z$ in $\Omega$, the classical theory of Cauchy-Lipschitz can be applied. 
Then there exists a maximal solution of the ordinary differential equation $V\in 
C^2([\dfrac{\delta}{2},A_\epsilon))$ such that 
$V(\dfrac{\tau}{2})=V_0(\dfrac{\tau}{2})$, $\dot{V}(\dfrac{\delta}{2})=\dot{V_0}(\dfrac{\delta}
{2})$ and $(V(x),\dot{V}(x))\in \Omega$ for all $x\in
 [\dfrac{\delta}{2},A_\epsilon)$. The local uniqueness given by the classical Cauchy-Lipschitz 
 theory ensures 
that we can connect $V_0$ and $V$ without 
 any problem and then get a solution of our ordinary differential equation on $[0;A_\epsilon)$. \hfill $\square$\\